%% file: ccacopf_conv.tex
\documentclass[journal]{IEEEtran}
\usepackage{cite}

\ifCLASSINFOpdf
   \usepackage[pdftex]{graphicx}
   \graphicspath{{./figs/}}
   \DeclareGraphicsExtensions{.pdf,.jpeg,.png}
\else
\fi
\usepackage{amsmath}
\usepackage{amsfonts}

\usepackage{algorithmic}
\usepackage{multirow}

\usepackage{xcolor}

\renewcommand{\b}[1]{\mathbf{#1}}
\newcommand{\bs}[1]{\boldsymbol{#1}}
\newcommand{\biidxex}[3]{\mathbf{#1}^{#3}_{#2}}
\newcommand{\iidxex}[3]{{#1}^{#3}_{#2}}
\newcommand{\bex}[2]{\mathbf{#1}^{#2}}
\newcommand{\biidx}[2]{\mathbf{#1}_{#2}}
\newcommand{\tp}{^{\top}}

\newcommand{\jbbo}[1]{\textcolor{black}{#1}} 
\newcommand{\jbb}[1]{\textcolor{black}{#1}} 

\newcommand{\jbbr}[1]{\textcolor{black}{#1}} 
\newcommand{\jbbR}[1]{\textcolor{black}{#1}} 

\begin{document}

\onecolumn

\input{cover}

\newpage

\twocolumn

%

\title{Convergence Analysis of Fixed Point Chance Constrained Optimal Power Flow Problems}
%
%
%

\author{Johannes J.~Brust, and~Mihai~Anitescu, \emph{Member}, \emph{IEEE}
\thanks{J.J. Brust and M. Anitescu were with Argonne National Laboratory, Mathematics and Computer Science Division
upon completion of this article, Lemont, IL, e-mails: jjbrust@ucsd.edu,anitescu@mcs.anl.gov}
\thanks{This work was supported by the U.S. Department of Energy, Office of Science,
Advanced Scientific Computing Research, under Contract DE-AC02-06CH11357
at Argonne National Laboratory.}} 

%
%

\ifCLASSOPTIONpeerreview
	\markboth{Draft I} 
	{Brust \MakeLowercase{\textit{et al.}}: Convergence of Chance Constrained Solver}
\fi
%



\maketitle

\begin{abstract}
For optimal power flow problems with chance constraints, a particularly effective
method is based on a fixed point iteration applied to a sequence of deterministic
power flow problems. However, a priori, the convergence of such an
approach is not necessarily guaranteed. This article analyses the convergence
conditions for this fixed point approach,
\jbbo{and reports numerical experiments including} \jbbo{for} \jbbo{large IEEE networks.}
\end{abstract}

\begin{IEEEkeywords}
Fixed point method, Chance Constraints, Stochastic Optimizations, AC Optimal Power Flow
\end{IEEEkeywords}

%
\IEEEpeerreviewmaketitle

\section{Introduction}
%
%
%
%
\IEEEPARstart{C}{hance} constrained optimization problems are
\jbbr{often} computationally very challenging. However, when modeling 
the effects of uncertainty on optimal power networks, potentially
large chance constrained optimization problems arise 
\cite{ZhangLi11,ZhangLi11a,HojjatJavidi15,BakerDallAneseSummers16,VrakopoulouEtal13}. 
\jbbr{Typically,} stochastic optimal power flow models are developed
by reformulating a widely accepted deterministic model.
One such ``classical" model is the AC optimal power flow model (AC-OPF) \cite{matpower11}. 
Because AC-OPF
has multiple degrees of freedom with respect to the problem variables,
various stochastic power flow paradigms can be \jbbr{derived} from it. In particular, 
one can define different subsets of variables as stochastic variables.
\jbbo{For instance, in \cite{LubinDvorkinRoald19} power generation is regarded 
as being stochastic (with constant demands), resulting in probabilistic objective functions.}
This article develops a chance constrained AC optimal power
flow model (CC-AC-OPF) in which the objective function
is deterministic, \jbbo{and the demands can be stochastic.}
This enables direct
interpretation of the meaning of the optimal objective function values, 
\jbbo{and can be more realistic, when in practice demand is more of
an uncertainty source as compared to supply.}
Yet, independent
of how the stochastic power flow model is formulated it 
typically yields a chance constrained optimization problem.
 In order to also solve large instances of the resulting 
CC-AC-OPF problems, a fixed point iteration
on a sequence of modified, related, and simpler deterministic AC-OPF problems
is solved. 
\jbbo{Iterative methods that solve a sequence of simpler problems have been very effective
on a variety of recent power systems problems \cite{bienstock2014chance,LeeTuritsynMolzahnRoald20,FrolovRoaldChertkov19}.
In the context of chance-constrained optimization, such an approach has been 
successfully used in \cite{ra18,SchmidliEtAl16,xuEtAl21} among more.}
However, prior to this work, there did not exist analytical criteria for when one
can expect this fixed point iteration to converge. Therefore, this article describes
a convergence analysis of the chance constrained fixed point iteration 
and tests the method on a set of standard IEEE networks. The numerical
experiments report results for networks with up to 9000 nodes.
\jbbo{To summarize, the contributions of this work are: (1) The formulation of a new
CC-AC-OPF model with a deterministic objective function, derived from uncertain demands; 
(2) The application and analysis of a fixed-point
algorithm for an implicit chance constrained problem. Even though
the use of iterative (fixed-point) techniques is widespread in power systems problems, previously no rigorous analysis had been undertaken
\jbbo{for this formulation}. 
In particular, we disentangle effects of model parameters and network properties on the convergence;
(3) We include numerical experiments on large test cases.}
\jbbr{The article is organized as follows: Sections I.A---D are preliminary and review the AC optimal power flow model
and how a chance-constrained model is obtained from it. Section II lists the fixed-point algorithm. We highlight that the reformulation and iterative solution of chance-constrained
AC-OPF has been proposed in \cite{ra18}, however a rigorous analysis of its convergence has not yet been developed. Therefore, Section III analyzes convergence properties of the fixed-point algorithm, while
numerical experiments are described in Section IV. Finally, we conclude with Section V.}




 

\subsection{AC Power Flow \jbbr{(Preliminary)}} 
The power network is defined by $ N $ nodes (buses) and $ N_{\text{line}} $
edges (lines). Associated to each bus is a voltage magnitude, $ v_i $,
frequency $ \theta_i $ and power
generation or consumption at the $i^{\text{th}}$ bus, for $ 1 \le i \le N $. In particular, let  $ p^{g}_i $
be the real power generated at bus $ i $, $ q^{g}_i $ the
corresponding reactive power and $p^{d}_i$, $q^{d}_i$ the real and reactive
demands, respectively. The buses are furthermore divided into two
groups: generators and loads. The indexing sets are 
$G$, and $L$, which are related to the set of all buses by $B = G \cup L$ (also $ N = N_G + N_L $).
Throughout, we will use the indexing sets $ \{G,L\} $ as subscripts to bold font vector variables. Note that in our setup each bus fits exactly into one of the two categories
(though one can easily set a virtual zero load at a node). 
It is customary to assume that the network contains one reference bus (typically at $i=1$).
 Load buses 
do not have active power generation and are thus defined by 
$ \biidxex{p}{L}{g} = \biidxex{q}{L}{g} = \b{0} $, where $ \bex{p}{g}, \bex{q}{g} $ are vectors
that contain the $p^{g}_i$'s and $q^{g}_i $'s. The power flow 
in the system is described by the power flow equations,
which couple the variables (e.g., \cite[Sections III-IV]{matpower11}). Specifically, let $ G_{ik}, B_{ik} $
denote the entries of the network's admittance matrix 
and define the quantities $ \theta^{ik} \equiv \theta_i - \theta_k $,
$ c_{ik} \equiv G_{ik} \cos(\theta^{ik}) + B_{ik}\sin(\theta^{ik}) $ and
$ d_{ik} \equiv G_{ik} \sin(\theta^{ik}) - B_{ik}\cos(\theta^{ik}) $. With 
these definitions, the $2N$ power flow equations, for $i = 1, \ldots, N$, are
\begin{align}
	\label{eq:pf}
	 & v_i \sum_{k=1}^N v_k c_{ik} - (\iidxex{p}{i}{g} - \iidxex{p}{i}{d}) = 0, \\ \nonumber
	 & v_i \sum_{k=1}^N v_k d_{ik} - (\iidxex{q}{i}{g} - \iidxex{q}{i}{d}) = 0,
\end{align}
\jbbo{when grouped by busses. }
If we let $ \b{d} = \text{vec}(\bex{p}{d},\bex{q}{d})  $ be the $ \mathbb{R}^{2N} $ vector containing
 $\iidxex{p}{i}{d},\iidxex{q}{i}{d}$ then in vector notation \eqref{eq:pf} is expressed as $\b{f}(\b{v},\bs{\theta},\bex{p}{g}, \bex{q}{g}; \b{d}) = \b{0} $,
 where $ \b{f} $ represents the nonlinear equations in \eqref{eq:pf}.
 Note, however, that $ \b{f} $ is linear in $ \b{d} $, a fact we will use later.
 \jbbo{Moreover, $ \b{d} $ is regarded as a parameter and not a variable.}


\subsection{AC Optimal Power Flow \jbbr{(Preliminary)}}
\label{sec:acopf}
Optimal power flow determines the best set of variables
that satisfy the network constraints. In addition to the 
power flow equations from \eqref{eq:pf}, branch current
bounds are typically included in the optimization problem. In particular,
let $LC$ be the set of all line connections, i.e., the set of index pairs that
describe all line connections (e.g., if bus $i$ is connected
to bus $k$, then $ (i,k) \in LC $). Subsequently, the so-called branch current constraints are 
$(D^{\text{max}}_{ik})^2 - (v_i \cos(\theta_i) - v_k \cos(\theta_k))^2 - (v_i \sin(\theta_i) - v_k \sin(\theta_k))^2 \ge 0$,
$\forall (i,k) \in LC $ for constant current limits $D^{\text{max}}_{ik}$. In vector notation these $N_{\text{line}}$ constraints are 
represented as $ \b{g}(\b{v},\bs{\theta}) \ge \b{0}  $. It is also
desirable to include ``hard" bounds on the generation variables,
such as $ \b{l}_p \le \biidxex{p}{G}{g} \le \b{u}_p  $, where $\b{l}_p$,
$\b{u}_p$ represent constant lower and upper bounds.
The AC optimal power flow (AC-OPF) problem, for a cost function $ C_0\jbbr{(\cdot)} $, is thus formulated
as
\begin{align}
	\label{eq:acopf}
	\underset{ \b{v}, \bs{\theta}, \bex{p}{g}, \bex{q}{g} }{ \text{ minimize } } 
	& C_0(\biidxex{p}{G}{g})\quad \text{subject to}  \\
	 \b{f}(\b{v},\bs{\theta},\bex{p}{g}, \bex{q}{g}; \b{d}) &= \b{0} \label{eq:pf_cons} \\
	 \b{g}(\b{v},\bs{\theta}) &\ge \b{0} \nonumber \\
	 \b{l}_p &\le \biidxex{p}{G}{g} \le \b{u}_p \nonumber \\
	 \b{l}_q &\le \biidxex{q}{G}{g} \le \b{u}_q \nonumber \\
	 \b{l}_{v} &\le \b{v} \le \b{u}_{v} \nonumber \\
	 \b{l}_{\theta} &\le \bs{\theta} \le \b{u}_{\theta} \nonumber	 	 
\end{align}
The cost function is typically a convex quadratic function that only depends
on the real power generation. Specifically, 
$ C_0(\biidxex{p}{G}{g}) = \sum_{i \in G} q_{ii} (\iidxex{p}{i}{g})^2 + \sum_{i \in G} q_{i} \iidxex{p}{i}{g} + q_{00}    $ for cost
data $ q_{ii}, q_i $ and $ q_{00} $. \jbbr{A local} solution $ \b{s}^* = \text{vec}(\b{v}^{*}, \bs{\theta}^{*}, (\bex{p}{g})^*, (\bex{q}{g})^*) $
of \eqref{eq:acopf} is called the optimal power flow point, or OPF point.

\subsection{Chance Constrained AC Optimal Power Flow \jbbr{(Preliminary)}}
\label{subsec:ccacopf}
In order to introduce uncertainty related to e.g., renewable energy
into the OPF problem \eqref{eq:acopf} we regard the demand terms in
\eqref{eq:pf} (i.e., $\iidxex{p}{i}{d}, \iidxex{q}{i}{d} $) as forecasted values with possible error. Specifically,
these stochastic quantities are represented as 
\begin{equation}	
	\label{eq:pd_rand}
	\iidxex{p}{i}{d} + \omega_i, \quad  \iidxex{q}{i}{d} + \omega_{N+i},
\end{equation}
where $\omega_i$ and $ \omega_{N+i} $ represent forecasting errors. For compactness,
the stochastic errors are represented by the $2N$ vector $ \bs{\omega} $. Here we assume
$ \bs{\omega} $ is normal, and relaxing this assumption yields different approaches. 
The chance constrained AC-OPF model in this Section
is developed such that the objective function only depends on deterministic variables. When the
objective function represents cost, deterministic values are meaningful and important. Since the 
power flow equations in \eqref{eq:pf} are overdetermined, i.e., $ \b{f}: \mathbb{R}^{2N+N_G} \to \mathbb{R}^{2N}  $,
this system has $N_G $ degrees of freedom. Subsequently, let $ \b{y} = \biidx{p}{G} $ 
represent a $ \mathbb{R}^{N_G} $ vector of deterministic variables and 
$ \b{x} = \b{x}(\bs{\omega}) = \text{vec}(\biidx{q}{G}(\bs{\omega}), \biidx{v}{L}(\bs{\omega}), \bs{\theta}(\bs{\omega})) $
a $ 2N $ vector of stochastic variables. The stochastic power flow equations are
represented as
\begin{equation}
	\label{eq:stochpfe}
	\b{f}(\b{x},\b{y}; \b{d}) + \bs{\omega} = \b{f}_{\omega}(\b{x}(\bs{\omega}),\b{y}; \b{d} + \bs{\omega}) \equiv \b{f}_{\omega} = \b{0} .
\end{equation}
If $ \bs{\omega} = \b{0} $ these equations reduce to the power flow equations in \eqref{eq:pf}.
Note that the power flow equations couple the variables and
that the uncertainty in $ \b{x} $ depends on the uncertainty in the demands ``$\b{d}$",
since $ \b{x} = \b{x}(\bs{\omega}) = \b{x}(\b{y}, \b{d} + \bs{\omega}) $. We set
$ C_0(\biidxex{p}{G}{g}) = C_0(\b{y}) $ to reflect the change in variables
for the stochastic optimal power flow problem and note that the objective function
is deterministic. Let $ \mathbb{P}( \b{z} \ge \b{0} ) \ge 1 - \bs{\epsilon} $ represent a vector of inequalities
with non-negative probability constraints, for which each element in the left hand side corresponds to a
cumulative probability (for $ z_i \ge 0 $) and each element of $ \bs{\epsilon} $ is in the interval $ 0 < \epsilon < 1 $. Then the
chance constrained (stochastic) optimal power flow problem is
given by
\begin{alignat}{2}
	\label{eq:ccacopf}
	\underset{ \b{x}(\bs{\omega}), \b{y} }{ \text{ minimize } } 
	& C_0(\b{y})\quad &\text{subject to}  \\
	 \b{f}_{\omega}(\b{x}(\bs{\omega}),\b{y}; \b{d} + \bs{\omega}) &= \b{0} \quad & \forall \bs{\omega}  \nonumber \\ 
	 \mathbb{P}(\b{g}(\b{x}(\bs{\omega}),\b{y}) \ge \b{0}) 	& \ge \b{1} - \bs{\epsilon}_g \label{eq:ccacopf_consg} \\
	 \mathbb{P}( \b{l}_x \le \b{x}(\bs{\omega}) \le \b{u}_x ) 	& \ge \b{1} - \bs{\epsilon}_x  \label{eq:ccacopf_consb} \\
	 \b{l}_{y} &\le \b{y} \le \b{u}_{y} \nonumber.	 	 
\end{alignat}
Observe that the problem in \eqref{eq:ccacopf} includes chance constraints (probability constraints)
on the line flow limits \eqref{eq:ccacopf_consg} and the stochastic variables $ \b{x}(\bs{\omega}) $ \eqref{eq:ccacopf_consb}, 
while deterministic limits are set on $ \b{y} $. Here $ \bs{\epsilon}_g $ and $ \bs{\epsilon}_x $ correspond to model
parameters for setting probability thresholds.

\subsubsection{Computing Chance Constraints}
To practically compute the chance constraints in problem \eqref{eq:ccacopf}, the stochastic variables are linearized around the 
error $ \bs{\omega} $ (zero mean). 
\jbbo{This is equivalent to assuming that $ \bs{\omega} $ is sufficiently small, which we proceed to do in the rest of the paper.}
In particular,
\begin{equation}
	\label{eq:linearize}
	\b{x}(\b{y},\b{d}+\bs{\omega}) \approx \b{x}(\b{y},\b{d}) + \frac{\partial \b{x}(\b{y},\b{d} + \bs{\omega} )}{\partial \bs{\omega}} \bigg|_{\bs{\omega = \b{0}}}  \bs{\omega} \equiv
	\biidx{x}{\omega}.
\end{equation}
Note that $ \b{x}(\b{y},\b{d}) = \b{x}_{0} $ and that $ \partial \b{x}(\b{y},\b{d}) \big / \partial \bs{\omega} $ are
deterministic. Thus the expectation
and variance of $ \biidx{x}{\omega} $ are $ \text{E}[\biidx{x}{\omega}] = \b{x}(\b{y},\b{d}) = \biidx{x}{0} $,
and $ \text{Var}[\biidx{x}{\omega}] = ( \partial \b{x}(\b{y},\b{d}) \big / \partial \bs{\omega} ) \text{Var}[\bs{\omega}] ( \partial \b{x}(\b{y},\b{d}) \big / \partial \bs{\omega} )^{\top} $, respectively. \jbbr{Alternatively, more accurate dynamics of the stochastic variables may be based on a $2^\text{nd}$ order expansion 
	$ \mathbf{x}(\b{y},\b{d} + \boldsymbol{\omega}) \approx \mathbf{x}(\b{y},\b{d}) + (\partial \mathbf{x} / \partial \boldsymbol{\omega})\boldsymbol{\omega} +
	\textnormal{``second order terms''}  $. When the covariance of the uncertainty has a particular structure (e.g., diagonal) then the mean of 
	the expansion can be determined analytically. The expansion's covariance is more involved and may need to be estimated.
	Another possibility to include nonlinearities might be a quadratic model of the load:
	$ \boldsymbol{\omega}_3 \texttt{.*} ( \mathbf{p}^d )^\texttt{.2} + \boldsymbol{\omega}_2 \texttt{.*} ( \mathbf{p}^d ) + \boldsymbol{\omega}_1 $, where 
	$ \texttt{.*} $ and $ ^\texttt{.2} $ are element-wise multiplications and squares}. 
For computational efficiency we use the probabilities of the linearized random variables. For instance, the constraints from \eqref{eq:ccacopf_consb} are written as
\begin{equation}
	\label{eq:cccons}
	\mathbb{P}(\biidx{x}{\omega} \le \b{u}_x) \ge \b{1} - \bs{\epsilon}_x, \quad \mathbb{P}(\b{l}_x \le \biidx{x}{\omega} ) \ge \b{1} - \bs{\epsilon}_x.
\end{equation}
To handle the vector of probabilities in \eqref{eq:cccons} one can use a Bonferroni bound \cite{mill81}, in which each
individual variable $ (\biidx{x}{\omega})_r $ for $ 1 \le r \le 2N $ satisfies an individual highly conservative probability constraint. 
However, when the variables are independent (or can be treated as such) then the probabilities can be separated without restrictions.
The mean
of $ (\biidx{x}{\omega})_r $ is $ (\biidx{x}{0})_r $, while the variance can also be explicitly computed. 
Define 
\begin{equation} 
	\label{eq:gamma}
	\partial \b{x} \big / \partial \bs{\omega} \equiv \bs{\Gamma},
\end{equation}
and let $ \b{e}_r $ be the r$^{\text{th}}$
column of the identity matrix. Moreover, denote
$ \text{Var}[\bs{\omega}] = \bs{\Sigma}^2 $. With this the variance of 
$ (\biidx{x}{\omega})_r $ is $ \left \| \b{e}^{\top}_r \bs{\Gamma} \bs{\Sigma} \right\|^2_2 $.
In turn, when the variables can be treated independently, individual probability constraints,
such as $ F^{\text{Nrm}}((\biidx{x}{\omega})_r \le (\b{u}_x)_r) \ge 1 - (\bs{\epsilon}_x)_r $ 
(where $ F^{\text{Nrm}} $ is the normal distribution function) can be represented as
\begin{equation*}
	((\b{u}_x)_r - (\biidx{x}{0})_r) \big / \left \| \b{e}^{\top}_r \bs{\Gamma} \bs{\Sigma} \right\|_2 \ge (F^{\jbbo{\text{Nrm}}})^{-1}(1- (\bs{\epsilon}_x)_r ), 
\end{equation*} 
where $ (F^{\jbbo{\text{Nrm}}})^{-1}( \cdot ) $ is the inverse cumulative distribution function.
Defining $ z_r = (F^{\jbbo{\text{Nrm}}})^{-1}(1- (\bs{\epsilon}_x)_r ) $ 
the constraints are represented as 
\begin{equation}
	\label{eq:lamcons}
	(\b{u}_x)_r - \lambda_r \ge (\biidx{x}{0})_r, \quad \lambda_r = z_r  \| \b{e}^{\top}_r \bs{\Gamma} \bs{\Sigma}  \|_2.
\end{equation}
Similarly, for $ 2N + 1 \le r_1 \le 2N + N_L $ \jbbo{and $ \bar{r}_1= r_1 - 2N $}, defining $ z_{r_1} = (F^{\jbbo{\text{Nrm}}})^{-1}(1- (\bs{\epsilon}_g)_{r_1}) $ the remaining 
probability constraints are represented as $ \b{g}(\biidx{x}{0},\b{y})_{\bar{r}_1} - \lambda_{r_1} \ge \b{0}$ with
\begin{equation}
	\label{eq:lamconsLF}
	\lambda_{r_1} = z_{r_1} \| \b{e}^{\top}_{\bar{r}_1} (\partial \b{g} \big / \partial \b{x}) \bs{\Gamma} \bs{\Sigma}  \|_2.
\end{equation}
Note that $ \lambda_r,\lambda_{r_1} $ depend on $ \b{x},\b{y} $ and $ \b{d} $, e.g., $ \lambda_r = \lambda_r(\b{x}(\b{y},\b{d})) $, 
which we will use later. Moreover the inequality in \eqref{eq:lamcons} is deterministic and thus
straightforward to compute once $ \bs{\Gamma} $ is known. Second, when $ (\bs{\epsilon}_x)_r $ is a small number 
(which is typically the case) then $ z_r > 0 $  and
$ \lambda_r > 0 $. Thus the $ \lambda_r,\lambda_{r_1} $ values are regarded as  
constraint tightenings and represent the effects of stochasticity in the constraints.
\jbbr{Note that when other distributions are desired, one can substitute the 
	$ (F^{\text{Nrm}})^{-1}(\cdot) $ c.d.f. (above \eqref{eq:lamcons} and elsewhere) with another one. Since
	investigations about distributional robustness have been previously conducted by other researchers we refer to 
	\cite[Sec. III.A]{Roald2017PowerSO} for 
	in-depth discussions.}

\subsubsection{Computing $ \partial \b{x} \big / \partial \bs{\omega} $}
\label{sec:derivsxo}
The partial derivatives $ \partial \b{x} \big / \partial \bs{\omega} = \bs{\Gamma} $ in \eqref{eq:lamcons}
are obtained by using the power flow equation
\begin{equation*}
	\frac{\partial }{ \bs{\partial \omega} } \left( \b{f}_{\omega}( \b{x}(\bs{\omega}), \b{y}; \b{d} + \bs{\omega} ) \right)  = 
	\frac{ \partial \b{f}_{\omega}}{\partial \b{x}} \frac{\partial \b{x}}{ \partial \bs{\omega}} + \frac{\partial \b{f}_{\omega}}{ \partial \bs{\omega}}= \b{0}.
\end{equation*}
First, note from \eqref{eq:stochpfe} that $ \partial \b{f}_{\omega} \big / \partial \bs{\omega} = \b{I} $. Second, the partial
derivatives are only needed at $ \bs{\omega} = \b{0} $, which yields the representation $ \partial \b{x} \big / \partial \bs{\omega} = - (\partial \b{f}_0 \big / \partial \b{x})^{-1}$ with the convention $ \b{f}_0 = \b{f} $. 
Since $ \b{x} = \text{vec}(\biidxex{q}{G}{g}, \biidx{v}{L}, \bs{\theta}) $ the so-called Jacobian matrix of
partial derivatives is given by
\begin{equation}
	\label{eq:jac}
	\frac{\partial \b{f}}{ \partial \b{x} } = 
		\left[
			\begin{array}{ c c c }
				\frac{\partial \b{f}}{ \partial \biidxex{q}{G}{g} } & \frac{\partial \b{f}}{ \partial \biidx{v}{L} } & \frac{\partial \b{f}}{ \partial \bs{\theta} }			
			\end{array}
		\right].
\end{equation}
The elements of this matrix can be computed from \eqref{eq:pf}. Note that only the last $ N $ equations in \eqref{eq:pf} depend on $ \bex{q}{g} $.
Subsequently, we define the indices $ 1 \le i \le N $ and $ j = N+ i $, 
as well as $ 1 \le g \le N_G $, $ 1 \le l \le N_L $ and $ 1 \le t \le N $. With this, the elements of 
the Jacobian from \eqref{eq:jac} are: 
\begin{align}
	\label{eq:jacelements}
	\left( \frac{\partial \b{f}}{ \partial \biidxex{q}{G}{g} } \right)_{i,g} &= 0,  \\
	\nonumber \left( \frac{\partial \b{f}}{ \partial \biidxex{q}{G}{g} } \right)_{j,g} &=
	\begin{cases}
		-1 & \text{ if } i = [G]_g \\
		0 & \text{ otherwise }
	\end{cases}	\\
	\nonumber
	\left( \frac{\partial \b{f}}{ \partial \biidx{v}{L} } \right)_{i,l} & = 
	\begin{cases}
		\sum_{k=1}^{N} v_{k} c_{ik} + 2c_{ii}v_i 	& \text{ if }  i = [L]_l  \\
		v_i c_{i[L]_l} 						& \text{ otherwise }
	\end{cases}  \\
	\nonumber \left( \frac{\partial \b{f}}{ \partial \biidx{v}{L} } \right)_{j,l} &=
	\begin{cases}
		\sum_{k=1}^{N} v_{k} d_{ik} + 2d_{ii}v_i 	& \text{ if }  i = [L]_l  \\
		v_i d_{i[L]_l} 						& \text{ otherwise }
	\end{cases}	\\
	\nonumber
	\left( \frac{\partial \b{f}}{ \partial \bs{\theta} } \right)_{i,t} & = 
	\begin{cases}
		-v_i\sum_{k=1}^{N} v_{k} d_{ik} 	~~~~~		& \text{ if }  i = t  \\
		v_i v_{t} d_{it} 	~~~~~		& \text{ otherwise }
	\end{cases}  \\
	\nonumber \left( \frac{\partial \b{f}}{ \partial \bs{\theta} } \right)_{j,t} &=
	\begin{cases}
		v_i\sum_{k=1}^{N} v_{k} c_{ik} 	~~~~~~~		& \text{ if }  i = t  \\
		-v_iv_{t} c_{it} 	~~~~~~~		& \text{ otherwise }
	\end{cases}
\end{align}
\jbb{The partial derivatives $ \partial \b{x} / \partial \bs{\omega} = ( \partial \b{f} / \partial \b{x} )^{-1} $ are defined by an inverse.
This inverse is typically well defined for regular optimal power flow problems, as described in \cite[Section III. B]{LubinDvorkinRoald19}
and \cite{Bienstock15} (\jbbr{if numerically the Jacobian matrix becomes (nearly) singular, it may be corrected by shifting its diagonal elements by adding a multiple of the identity matrix}). Therefore, the smallest singular value of $ \partial \b{f} / \partial \b{x} \equiv \b{J} $ is positive, i.e.,
$ \sigma_{\text{min}}( \b{J} ) > 0 $. A positive lower bound for the smallest singular value of a matrix is described in \cite[Theorem 1]{HongPan92}.
Let $ \b{J}_{:,i} $ denote the $i^{\text{th}}$ column of $ \b{J} $ and let $ \b{J}_{i,:} $ be the $i^{\text{th}}$ row of $ \b{J} $.
Then a lower bound for the smallest singular value is:
\begin{equation*}
\sigma_{\text{min}}(\b{J}) \ge  \hat{K}_{\Gamma} > 0,
\end{equation*}
where
$ \hat{K}_{\Gamma} = \left( \frac{\hat{n}-1}{\hat{n}} \right)^{\frac{\hat{n}-1}{2}} |\b{J} | 
\text{ max}\left(\frac{\text{min}(\|\b{J}_{:,i} \|_2) }{ \prod_i \| \b{J}_{:,i} \|_2 }, \frac{\text{min}(\|\b{J}_{i,:} \|_2) }{ \prod_i \| \b{J}_{i,:} \|_2 } \right),
$
$ \hat{n} \equiv 2N $, and where the determinant is $ | \b{J} | = \text{det}(\b{J}) $. }
This means that
\begin{equation}
	\label{eq:jacbound}
		\| (\partial \b{f} \big / \partial \b{x})^{-1} \|_2 
		= \jbb{ \| \b{J}^{-1} \|_2  = \frac{1}{ \sigma_{\text{min}} (\b{J} )  } \le 1 \big /  \hat{K}_{\Gamma} \equiv K_{\Gamma} }
\end{equation}
where $ \hat{K}_{\Gamma}, K_{\Gamma} $ are finite constants. \jbb{In order to compute solves
with $ \b{J} $ (which is part of computing the constraints in \eqref{eq:lamcons} and \eqref{eq:lamconsLF}) the LU factorization is based on the decomposition $ \b{J} = \b{P} \b{L} \b{U} $,
where $ \b{P} $ is a permutation matrix, $ \b{L}$ is a unit lower triangular matrix
and $ \b{U} $ is an upper triangular matrix. The determinant and the bounds in \eqref{eq:jacbound} are thus available ``without extra expense'' based on solves with $ \b{J} $, by multiplying the
diagonal elements of $ \b{U} $, since $ | \b{J} | = | \b{U} |  $ and $ \b{U} $ is upper triangular.}
\jbbr{Since the determinant can often become large (even if a matrix is well conditioned) a possibly
preferable approach of computing constant $ K_{\Gamma} $ is to exploit the inequality $ \| \b{J}^{-1} \| \le \sqrt{ \|\b{J}^{-1} \|_1 \|\b{J}^{-1} \|_{\infty}  } \equiv K_{\Gamma} $.
Note that computing $ K_{\Gamma} $ can be inexpensive since $ \b{J}^{-1} $ is computed as part of e.g., \eqref{eq:lamcons}.}

\subsection{Implicit Chance Constrained Optimal Power Flow}
\label{subsec:iccopf}
An optimal power flow problem
that combines components of the AC-OPF problem in \eqref{eq:acopf} and 
of the chance constrained problem in \eqref{eq:ccacopf} is obtained by using the constraints
from \eqref{eq:lamcons} and \eqref{eq:lamconsLF}. This reformulated problem
incorporates stochastic effects, while at the same time enables efficient computations.
The corresponding \emph{implicit} chance-constrained AC-OPF problem is:
\begin{align}
	\label{eq:iacopf}
	\underset{ \b{v}, \bs{\theta}, \bex{p}{g}, \bex{q}{g} }{ \text{ minimize } } 
	& C_0(\biidxex{p}{G}{g})\quad \text{subject to}  \\
	 \b{f}(\b{v},\bs{\theta},\bex{p}{g}, \bex{q}{g}; \b{d}) &= \b{0} \nonumber \\
	 \b{g}(\b{v},\bs{\theta}) &\ge \bs{\lambda}_{g}(\biidxex{p}{G}{g}) \nonumber \\	 
	 \b{l}_q + \bs{\lambda}_q(\biidxex{p}{G}{g}) &\le \biidxex{q}{G}{g} \le \b{u}_q - \bs{\lambda}_q(\biidxex{p}{G}{g}) \nonumber \\
	 \b{l}_{v} + \bs{\lambda}_v(\biidxex{p}{G}{g}) &\le \biidx{v}{L} \le \b{u}_{v} - \bs{\lambda}_v(\biidxex{p}{G}{g}) \nonumber \\
	 \b{l}_{\theta} + \bs{\lambda}_{\theta}(\biidxex{p}{G}{g}) &\le \bs{\theta} \le \b{u}_{\theta} - \bs{\lambda}_{\theta}(\biidxex{p}{G}{g}) \nonumber \\
	 \b{l}_p &\le \biidxex{p}{G}{g} \le \b{u}_p \nonumber,
\end{align}
where  $ \bs{\lambda}_q(\biidxex{p}{G}{g}),\bs{\lambda}_v(\biidxex{p}{G}{g}),\bs{\lambda}_{\theta}(\biidxex{p}{G}{g}) $ are computed
using \eqref{eq:lamcons} and $ \bs{\lambda}_g(\biidxex{p}{G}{g}) $ is computed using \eqref{eq:lamconsLF}. 
The problem can be seen to be implicit with regards to probability constraints, because the effects of uncertainty are implicitly included in the tightenings of
some constraints by the non-negative $ \bs{\lambda}$'s.

\section{Method}
\label{sec:method}
The solution of the potentially large nonlinear optimization problem in \eqref{eq:iacopf} can be computed
directly. However, the computation of the $ \bs{\lambda}$'s adds nonlinearities, because they depend on the
matrix $ \bs{\Gamma} = - (\partial \b{x} / \partial \bs{\omega})^{-1} \in \mathbb{R}^{2N \times 2N} $, which is an inverse. Because of this,
the problem in \eqref{eq:iacopf} is still computationally difficult. Instead of solving
\eqref{eq:iacopf} directly,
\jbbr{an iterative scheme, which computes an approximate solution of \eqref{eq:iacopf}, 
by solving a sequence of simpler problems has been proposed in  \cite{ra18}}.
Such an algorithm is stated as:

\vspace{0.2cm}
\fbox{\parbox{0.90\linewidth}{
 \textbf{Algorithm 1}
\label{alg:alg1}
\begin{algorithmic}[1]
\STATE{Inputs: $ \b{s}^{(0)} = \text{vec}(\b{v}^{(0)},\bs{\theta}^{(0)}, \b{p}^{(0)}, \b{q}^{(0)} )$, $0 < \tau_q$, $0 < \tau_v$,
$0 < \tau_{\theta}$, $0 < \tau_{g}$ 
$0 < \epsilon_q \le 0.5$, $0 < \epsilon_v \le 0.5$,
$0 < \epsilon_{\theta} \le 0.5$, $0 < \epsilon_{g} \le 0.5$, $ 0 \prec \bs{\Sigma} $, $\bs{\lambda}^{(0)}_q = \bs{\lambda}^{(0)}_v = \bs{\lambda}^{(0)}_{\theta} = \bs{\lambda}^{(0)}_{g} = 0, k = 0$;}
\STATE{for $ k = 0,1,\ldots $} 
\STATE{\label{alg:alg1_s}Solve \eqref{eq:iacopf} with $\bs{\lambda}^{(k)}_q$,$\bs{\lambda}^{(k)}_v$,$\bs{\lambda}^{(k)}_{\theta}$ fixed and obtain $ \b{s}^{(k+1)} $;}
\STATE{\label{alg:alg1_lam} Compute $\bs{\lambda}^{(k+1)}_q$,$\bs{\lambda}^{(k+1)}_v$,$\bs{\lambda}^{(k+1)}_{\theta}$, $\bs{\lambda}^{(k+1)}_{g}$ using \eqref{eq:lamcons} and \eqref{eq:lamconsLF};} 
\IF{\begin{align*} 
\| \bs{\lambda}^{(k+1)}_q - \bs{\lambda}^{(k)}_q \|_{\infty} &\le \tau_q  \text{ and } \\
\| \bs{\lambda}^{(k+1)}_v - \bs{\lambda}^{(k)}_v \|_{\infty}  &\le \tau_v  \text{ and } \\
\| \bs{\lambda}^{(k+1)}_{\theta} - \bs{\lambda}^{(k)}_{\theta} \|_{\infty} &\le \tau_{\theta} \text{ and } \\
\| \bs{\lambda}^{(k+1)}_{g} - \bs{\lambda}^{(k)}_{g} \|_{\infty} &\le \tau_{g}
\end{align*}}
\STATE{Stop. Return $\b{s}^{(k+1)}$;}
\ELSE
\STATE{$ k = k +1$;}
\ENDIF
\end{algorithmic}}}

\vspace{0.2cm}

Note that Algorithm 1 stops when the changes in the $ \bs{\lambda}$'s become small. However
no criteria have yet been specified for when one can expect the iteration to converge.
Therefore, we analyze conditions for which one
can expect Algorithm 1 to converge. 

\section{Analysis}
\label{sec:analysis}

The analysis of Algorithm 1 is based on the \jbbr{insight} that this
iterative algorithm is a fixed point iteration. Therefore, in order to derive its convergence conditions, we investigate
the convergence conditions of the fixed point iteration. Throughout this section, the problem from \eqref{eq:iacopf}
is reformulated as 
\begin{align}
	\label{eq:iacopf_eq}
	& \underset{ \b{s} }{\text{minimize } } 
	C_1( \b{s} )\quad \text{subject to}  \\
	& \jbbo{\b{f}(\b{s}) = \b{0}}, \quad \jbbo{\b{h}(\b{s};\bs{\lambda})} \jbbo{\ge \b{0}} \nonumber
\end{align}
with $ \b{s} = \text{vec}(\b{v}, \bs{\theta}, \bex{p}{g},\bex{q}{g}) $, $ C_1(\b{s}) = C_0( \b{p}^g_G )$,
$ \bs{\lambda} = \bs{\lambda}(\b{s}) =  \text{vec}(\bs{\lambda}_q,\bs{\lambda}_v,\bs{\lambda}_{\theta},\bs{\lambda}_g)$, 
$ \jbbo{\b{h}(\b{s};\bs{\lambda}) \ge \b{0} } \in \mathbb{R}^{m} $ ($m = N_L + 4N$ inequality constraints in \eqref{eq:iacopf}) and $ \b{f}(\b{s};\b{d}) = \b{f}(\b{s}) $.
\jbbo{Our analysis is consistent with classical nonlinear programming sensitivity results \cite{FiaccoKyparisis85,Fiacco80},
however is different because of the fixed point component of the algorithm. In classical nonlinear programming
$ \bs{\lambda} $ is taken as an independent perturbation parameter, with focus on the sensitivities of $ \b{s} $
with respect to $ \bs{\lambda} $. 
However, because the (implicit) chance-constrained problem also contains the dependence $ \bs{\lambda} = \bs{\lambda}(\b{s}) $
we are led to also} \jbbo{acknowledge the} \jbbo{interplay of $ \bs{\lambda} $ and $\b{s}$.}
Recall that a fixed point, say $ \bs{\lambda}^* $, of a continuous function/mapping, say $ \jbbr{M}(\bs{\lambda}) $, can be defined by the two conditions
\begin{equation}
\label{eq:convcond}
\underset{\text{Condition 1}}{\jbbr{M}(\bs{\lambda}^*) = \bs{\lambda}^*},  \quad 
\underset{\text{Condition 2} }{\| \partial \jbbr{M}(\bs{\lambda}^*) \big / \partial \bs{\lambda}  \|_2 \in [0,1)}
\end{equation}
\jbbr{(see for instance \cite[Sec. 10.1]{burden2010numerical} and \cite[Sec. 6.3]{strogatz} on fixed points.)} \jbbr{These two conditions are generalizations of a scalar fixed-point (FP) to vector valued functions. Condition 1 represents the definition of a FP, which is a point that remains unchanged after being mapped. Condition 2 is 
a sufficient condition for a fixed point mapping to converge \cite[Theorem 10.6]{burden2010numerical}). This condition implies that recursive mappings of a fixed-point are characterized by vanishing 
derivatives (successively applying the mapping results in a stable point, which make such mapping be sometimes called a contraction \cite[Theorem 10.6]{burden2010numerical}). Note that Algorithm 1 in line 5 checks whether the differences in successive constraint tightenings are small, i.e., this is a numerical check of Condition 1. 
\jbbR{We further analyze properties of Condition 2 to obtain insights into the convergence behavior of the algorithm.}
} 
Now, to closer investigate the relation between the variables in Algorithm 1, define the mapping
in line \ref{alg:alg1_s} that determines $ \b{s}^{(k+1)} $ from $ \bs{\lambda}^{(k)} $ by $ \b{s}_M(\bs{\lambda}^{(k)}) $ (i.e., $ \b{s}^{(k+1)} = \b{s}_M(\bs{\lambda}^{(k)}) $), and
the operation in line \ref{alg:alg1_lam} that determines $ \bs{\lambda}^{(k+1)} $ from $ \b{s}^{(k+1)} $ by $ \bs{\lambda}_M(\b{s}^{(k+1)}) $
(i.e., $ \bs{\lambda}^{(k+1)} = \bs{\lambda}_M(\b{s}^{(k+1)}) $).
Then the statements in the algorithm recursively define the next iterates, starting from $ k =0 $, $ \b{s}^{(k)}$, $ \bs{\lambda}^{(k)}$,
so that, 
\begin{alignat*}{3}
	\b{s}^{(k+1)} 	&= \b{s}_{M}( \bs{\lambda}^{(k)} ), \quad  			&& \bs{\lambda}^{(k+1)} 	&&= \bs{\lambda}_{M}( \b{s}^{(k+1)} )\\
				&= \b{s}_{M}(\bs{\lambda}_{M}( \b{s}^{(k)} )), \quad 	&&			&&= \bs{\lambda}_{M}( \b{s}_{M}( \bs{\lambda}^{(k)} )).
\end{alignat*}
This means that there is a mapping that generates $ \b{s}^{(k+1)} $ from $ \b{s}^{(k)} $ and one that
generates $ \bs{\lambda}^{(k+1)} $ from $ \bs{\lambda}^{(k)} $. In particular, if $ \bs{\lambda}^{(k)} \to \bs{\lambda}^* $ then $ \b{s}^{(k)} \to \b{s}^* $
and vice-versa. \jbbr{Our analysis is based on the assertion that if Condition 2 in \eqref{eq:convcond} holds, then $ \bs{\lambda}^* $ (and $ \b{s}^* $) is a fixed point.}

\jbbr{First, we describe the basic properties of a solution to \eqref{eq:iacopf_eq}, which enables the use of results
from nonlinear programming later on, and also the establishing of reasonable assumptions.}


\subsection{\jbbr{Basic Conditions}}
\label{sec:basic}
\jbbr{A} solution to the nonlinear programming problem \eqref{eq:iacopf_eq} is characterized by a set of
conditions \jbbr{for the} Lagrangian:
\begin{equation*}
	L(\b{s},\bs{\mu},\bs{\rho};\bs{\lambda}) =
	C_{\jbbr{1}}(\b{s}) + \bs{\mu}^\top \b{f}(\b{s}) + \bs{\rho}^\top \b{h}(\b{s};\bs{\lambda}),
\end{equation*}
where  $ \bs{\mu} \in \mathbb{R}^{2N} $ and $ \bs{\rho} \ge \b{0} \in \mathbb{R}^m $.
Subsequently, the Karush-Kuhn-Tucker (KKT) \cite{Karu39,KuhnTucker51} optimality conditions
are the set of nonlinear conditions that define a solution of \eqref{eq:iacopf_eq}:
\begin{align}
	\label{eq:KKT}
	\frac{\partial}{\partial \b{s}} L(\b{s},\bs{\mu},\bs{\rho};\bs{\lambda}) &= \b{0} \nonumber \\
	\b{f}(\b{s}) &= \b{0} \nonumber \\
	\b{h}(\b{s},\bs{\lambda}) &\ge \b{0} \\
	\bs{\rho}_i \b{h}(\b{s},\bs{\lambda})_i &= 0, \quad i =1,\cdots,m \nonumber \\
	\bs{\rho} &\ge \b{0}. \nonumber
\end{align}
The set of active inequality constraints is defined as
\begin{equation*}
	\mathcal{A}(\b{s}) = \{ i: \b{h}(\b{s},\bs{\lambda})_i = 0 \}.
\end{equation*}
\jbbr{A solution to \eqref{eq:KKT} can be found when the columns}
of the constraint Jacobians are linearly independent, i.e., when 
\begin{equation}
	\label{eq:linIndep}
	\left[ \quad \frac{\partial}{ \partial \b{s} } \b{f}(\b{s}), \quad  \frac{\partial}{ \partial \b{s} } \b{h}(\b{s},\bs{\lambda})_i \quad \right], \quad i \in \mathcal{A}(\bs{s}),
\end{equation}
are linearly independent. Moreover, second order conditions (which state
that the Lagrangian Hessian is positive definite in the nullspace of the constraint derivatives) 
ensure strict local optimality of a KKT point. Finally, if the active set is unchanged in a small neighborhood
around $ \bs{\lambda} $, then changes in $ \b{s} $ with regards to changes in $ \bs{\lambda} $ are continuous. Summarizing
we assume the three conditions.
\jbbr{Assumptions (Analysis):}
\begin{flushleft}
\begin{tabular}{l l}
	\jbbr{A.1}:&Linear independence in the constraint Jacobians \eqref{eq:linIndep} \\
	\jbbr{A.2}:&Second order sufficient conditions hold for \eqref{eq:iacopf_eq} \\
	\jbbr{A.3}:&Strict complementary slackness:
	$ \bs{\rho}_i > 0, i \in \mathcal{A}(\b{s}) $.
\end{tabular}
\end{flushleft}

The first two assumptions state that a local minimum for the problem in \eqref{eq:iacopf_eq} exists.
The third ensures continuity of derivatives, such as  $ \frac{\partial \b{s}}{ \partial \bs{\lambda}} $
and $ \frac{\bs{\lambda}_M}{ \partial \bs{\lambda} } $.
\jbbr{We further assume access to a solver that can find a local minimum when it exists. Because of the
continuity of partial derivatives the active set at a (local) solution does not change. A numerical investigation of a zig-zag behavior
of Algorithm 1 due to a disconnected feasible set can be found in \cite{Roald2017PowerSO}.
Note that when changes in the constraint tightenings are not abrupt, but smooth, assuming continuity is reasonable. 
We like to note that \cite[Sec. V.A]{Roald2017PowerSO} empirically observed that in cases where tightenings are small (relative to other terms)
the algorithm converged more frequently on linear programs.} 

\subsection{\jbbR{Analysis Outline}}
\label{sec:anoutline}
\jbbR{A goal of the analysis is to deduce properties of the quantity $ \| \frac{\partial \bs{\lambda}_M}{ \partial \bs{\lambda} } \| $ (in order to investigate when
Condition 2 in \eqref{eq:convcond} will hold). To achieve this first a representation of $ \frac{\partial \bs{\lambda}_M}{ \partial \bs{\lambda} } $ is needed, which is derived in Part I. 
Subsequently, the representation of $ \frac{\partial \bs{\lambda}_M}{ \partial \bs{\lambda} } $ depends on changes in variables and constraints that are developed
as the sensitivities $ \frac{\partial \b{s}}{ \partial \lambda_n} $ and $ \frac{\partial^2 \b{f}}{ \partial \lambda_n \partial \b{x}  } $ in Part II. Ultimately,
Part III uses the previous outcomes to deduce $ \left\| \frac{\partial \bs{\lambda}_M}{ \partial \bs{\lambda} } \right \|_2 $  and implications.}
(We use the notation $ (\bs{\lambda})_n = \lambda_n $ to denote the $n^{\text{th}}$ element of $ \bs{\lambda} $)
\begin{flushleft}
\begin{tabular}{l l}
	Part I: &  Representation of $ \frac{\partial \bs{\lambda}_M}{ \partial \bs{\lambda} }  $ \\ 
	Part II: & Sensitivities $ \frac{\partial \b{s}}{ \partial \lambda_n} $ and $ \frac{\partial^2 \b{f}}{ \partial \lambda_n \partial \b{x}  } $ \\ 
	Part III: & Bound on $ \left\| \frac{\partial \bs{\lambda}_M}{ \partial \bs{\lambda} } \right \|_2 $  and Implications 
\end{tabular}
\end{flushleft}
\jbbr{To develop a representation for $ \frac{\partial \bs{\lambda}_M}{\partial \bs{\lambda}} $ }
we analyze the expression
\begin{equation*}
	(\bs{\lambda}_M(\b{s}(\bs{\lambda})))_r = z_r \| \b{e}^{\top}_r \bs{\Gamma}(\b{s}(\bs{\lambda})) \bs{\Sigma}  \|_2,
\end{equation*}
which is \eqref{eq:lamcons} with explicit dependencies on $ \b{s} $ and $ \bs{\lambda} $. (The conditions for \eqref{eq:lamconsLF} are
done in a similar way, with an additional constant for the derivatives of the line flow constraints: $ \| \partial \b{g}  \big/ \partial \b{x} \|_2 \le K_g $).

\subsection{Part I: Representation of $ \frac{\partial \bs{\lambda}_M}{ \partial \bs{\lambda} }  $}
\label{sec:part_I}
\jbbr{In this section we first make the relation between $ \b{s} $ and $ \bs{\lambda} $ explicit. Thus}
the $ \bs{\Gamma} $ matrix is represented as a function of $\b{s}$ and $ \bs{\lambda} $  i.e.,
		$
		\bs{\Gamma}(\b{s}(\bs{\lambda})) = - \b{J}(\b{s}(\bs{\lambda}))^{-1},
		$	
		where $ \b{J}(\b{s}(\bs{\lambda})) = \b{J} = \partial \b{f} \big / \partial \b{x} $ (cf. \eqref{eq:jac}).
%
For notation, we will be using the following definition
		\begin{equation}
		\label{eq:si}
		\b{w}_r (\b{s}(\bs{\lambda})) \equiv  \bs{\Sigma} \bs{\Gamma}(\b{s}(\bs{\lambda}))^{\top} \b{e}_r,
		\end{equation}
		and write
		$
		(\bs{\lambda}_M(\b{s}(\bs{\lambda})))_r = z_r[\b{w}_r (\b{s}(\bs{\lambda}))^\top \b{w}_r (\b{s}(\bs{\lambda}))]^{1/2}.
		$
		Because Condition 2 in \eqref{eq:convcond} involves partial derivatives, we take the partial derivative
		with respect to $ \lambda_n $ (the `$n$th' tightening), so that
		$
		\frac{\partial }{\partial \lambda_n} (\bs{\lambda}_M(\b{s}(\bs{\lambda})))_r = \frac{z_r^2}{(\lambda_M)_r}\left[ \b{w}_r^{\top}\frac{\partial \b{w}_r}{\partial \lambda_n}\right], 
		$
		denoting $ (\bs{\lambda}_M(\b{s}(\bs{\lambda})))_r = (\lambda_M)_r $ in the right hand side.
		This expression provides a basis for what the matrix of partial derivatives will look like. 
		For indices $ 1 \le r,n \le 2N $ the elements of the matrix of \emph{sensitivities} 
		with respect to the vector $ \bs{\lambda} $ is
				\begin{equation}
			\label{Jacob}	
			{\small}
			\left(\frac{\partial \bs{\lambda}_M}{\partial \bs{\lambda}}\right)_{rn} = 
			\frac{z_r^2}{(\lambda_M)_r}\left[ \b{w}_r^{\top}\frac{\partial \b{w}_r}{\partial \lambda_n}\right].
		\end{equation}	
		Note that since $\b{w}_r^{\top}$ is a row vector and $\frac{\partial \b{w}_r}{\partial \lambda_n}$ is a column vector the entries 
		in \eqref{Jacob} can be written as
		\begin{equation}
			\label{dotproduct} 
			\b{w}_r^{\top}\frac{\partial \b{w}_r}{\partial \lambda_n} =  \left\| \b{w}_r\right\|_2 \left\| \frac{\partial \b{w}_r}{\partial \lambda_n}\right \|_2\cos(\phi_{rn}),
		\end{equation}
		where $ \phi_{rn} $ represents the angle between $ \b{w}_r $ and $ \frac{\partial \b{w}_r}{\partial \lambda_n} $.
		Since $(\lambda_M)_r = z_r \left\| \b{w}_r \right\|_2$ 
		therefore the elements of the matrix in \eqref{Jacob} are
		$$
		\jbbr{
		 \frac{z^2_r}{(\lambda_M)_r} \left[ \b{w}_r^{\top}\frac{\partial \b{w}_r}{\partial \lambda_n} \right] 
		 = z_r  \left\| \frac{\partial \b{w}_r}{\partial \lambda_n}\right\|_2 \cos(\phi_{rn}).} 
		 $$
		 From \eqref{eq:si} it holds that 
		 $ 
		 \left\| \frac{\partial \b{w}_r}{\partial \lambda_n}\right\|_2 =  \left\| \bs{\Sigma} \frac{\partial \bs{\Gamma}^{\top}}{\partial \lambda_n} 			\b{e}_r\right\|_2 
		 $,		
which yields the subsequent inequalities
\begin{equation}
	\label{eq:boundjacel}
			\left(\frac{\partial \bs{\lambda}_M}{\partial \bs{\lambda}} \right)_{rn} \leq | z_r | \left\| \bs{\Sigma} \right\|_2 \left\|\frac{\partial \bs{\Gamma}\tp}{\partial \lambda_n} \b{e}_r \right\|_2. 
\end{equation}
Note that $ z_r $ and $ \bs{\Sigma} $ are constants, with upper bounds 
$ z_r \le \underset{ 1 \le r \le 2N }{\max | z_r |} \equiv K_1 $ and $ \| \bs{\Sigma} \|_2 \le K_2 $. Therefore,
a bound for $ \left\| \frac{\partial \bs{\Gamma}\tp}{\partial \lambda_n} \b{e}_r \right\|_2 $ fully specifies \eqref{eq:boundjacel},
and thus the elements of $ \frac{\partial \bs{\lambda}_M}{\partial \bs{\lambda}} $. 
\subsection{Part II: Sensitivities $ \frac{\partial \b{s}}{ \partial \lambda_n} $ and $ \frac{\partial^2 \b{f}}{ \partial \lambda_n \partial \b{x}  } $}
\label{subsec:part_II}
In order to derive a bound on $ \| \frac{ \partial \bs{\Gamma} \tp}{ \partial \lambda_n } \b{e}_r  \|_2 = \| \b{e}\tp_r \frac{ \partial \bs{\Gamma} }{ \partial \lambda_n } \|_2 $, recall that 
$ \jbbr{\bs{\Gamma} = \b{J}^{-1}}$. 
\jbbr{Therefore}, note that we can make use of the identity 
\begin{equation*}
	\frac{ \partial }{ \partial \lambda_n } \b{J}^{-1} = - \b{J}^{-1} \left(\frac{ \partial  }{ \partial \lambda_n } \b{J} \right)\b{J}^{-1}.
\end{equation*}
Specifically, defining the vector $ \b{j}\tp_r \equiv \b{e}\tp_r \b{J}^{-1} $ then
\begin{equation}
	\label{eq:gammaDerivBound}
	\left\| \b{e}\tp_r \frac{\partial  \b{J}^{-1}}{ \partial \lambda_n } \right \|_2
	= \left \| \b{j}\tp_r \frac{\partial \b{J}}{ \partial \lambda_n }  \b{J}^{-1} \right\|_2
	\le \left \| \b{J}^{-1} \right \|_2 \left \| \b{j}\tp_r \frac{\partial \b{J}}{ \partial \lambda_n }  \right\|_2
\end{equation}
With a bound on $ \left\| \jbbr{\b{J}}^{-1} \right\|_2 \le K_{\Gamma}$ (from \eqref{eq:jacbound}), it remains to derive an upper bound on 
$ \left\|  \b{j}\tp_r  \jbbr{\frac{\partial \b{J}}{\partial \lambda_n}} \right\|_2 $. 
Next we describe the essential components for such a bound, based on
the elements of the $2^{\text{nd}}$ derivative matrices
\begin{equation*}
	\frac{\partial }{\partial \lambda_n} \b{J}(\b{s}(\bs{\lambda})) = \frac{\partial }{\partial \lambda_n} \left( \frac{\partial \b{f}}{ \partial \b{x} } \right) = 
		\left[
			\begin{array}{ c c c }
				\frac{\partial^2 \b{f}}{ \partial \lambda_n \partial \biidxex{q}{G}{g} } & \frac{\partial^2 \b{f}}{ \partial \lambda_n \partial \biidx{v}{L} } & \frac{\partial^2 \b{f}}{ \partial \lambda_n \partial \bs{\theta} }			
			\end{array}
		\right].
\end{equation*}
The elements of $ \frac{\partial }{\partial \lambda_n} \left( \frac{\partial \b{f}}{ \partial \b{x} } \right) $ can be
computed from \eqref{eq:jacelements} once 
\begin{equation}
	\label{eq:lamsensitivities}
	\frac{\partial}{ \partial \lambda_n }v_i \equiv \partial_n v_i, \quad \jbbr{\text{and}} \quad \frac{\partial}{ \partial \lambda_n  }\theta_i \equiv \partial_n \theta_i,
\end{equation}
for $ 1 \le n \le 2N $, $ 1 \le i \le N $ are known. This is the case, because $ \frac{\partial}{ \partial \lambda_n } (v_i c_{ki}) $ and
$ \frac{\partial}{ \partial \lambda_n } (v_i d_{ki}) $, for $ 1 \le k \le N  $, (in e.g., \eqref{eq:jacelements}) can be computed from
these quantities. \jbbr{In order to develop expressions for the derivatives in \eqref{eq:lamsensitivities} we use the vector variable $ \b{s} $ (from \eqref{eq:iacopf_eq}), which contain elements $ v_i $ and $ \theta_i $. \jbbr{Subsequently results} from classical nonlinear programming theory,
\jbbr{based on assumptions A.1 --- A.3, ensure}
 the existence of a parametrized solution to \eqref{eq:iacopf_eq} with dependence on $ \bs{\lambda} $,
defined by}
\begin{equation*}
	\jbbr{
	\b{a}(\bs{\lambda}) \equiv
	\left[\:\: \b{s}(\bs{\lambda})\tp \:\: \bs{\mu}(\bs{\lambda})\tp \:\: \bs{\rho}(\bs{\lambda})\tp \:\: \right]\tp.}
\end{equation*}
Selecting a set of rows in $ \b{a}(\bs{\lambda}) $ extracts $ \b{s}(\bs{\lambda}) $ \jbbr{and thus also $ v_i $ and $ \theta_i $.}
The derivative of $ \b{a}(\bs{\lambda}) $ w.r.t. $ \bs{\lambda} $ can be computed from the KKT conditions
\begin{equation}
	\label{eq:KKTlam}	
	\left[
		\begin{array}{c}
			\frac{\partial}{ \partial \b{s} } L(\b{a}(\bs{\lambda}), \bs{\lambda}) \\ 
			\b{f}(\b{s}(\bs{\lambda})) \\
			\bs{\rho}_i \b{h}(\b{s}(\bs{\lambda}),\bs{\lambda})_i
		\end{array}
	\right]
	 \equiv \b{G}(\b{a}(\bs{\lambda}),\bs{\lambda})
	=
	 \b{0}, \quad i \in \mathcal{A}(\b{s}).
\end{equation}
From \eqref{eq:KKTlam} it holds that
\begin{equation*}
	\frac{\partial}{ \partial \lambda_n } \b{G}(\b{a}(\bs{\lambda}),\bs{\lambda}) =
	\frac{\partial \b{G}}{ \partial \b{a} } \frac{\partial \b{a}}{ \partial \lambda_n } + 
	\frac{\partial \b{G}}{ \partial \lambda_n } = \b{0},
	\end{equation*}
and $ \frac{\partial\b{a}}{\partial \lambda_n} = - \big[ \frac{\partial \b{G}}{\partial \b{a}} \big]^{-1} \frac{\partial \b{G}}{ \partial \lambda_n } $. 
Thus $ \frac{\partial \b{s}}{\partial \lambda_n} $ is obtained by extracting elements from $ \frac{\partial \b{a}}{ \partial \lambda_n } $.
The existence of the inverse and hence the derivatives of $ \b{a} $ are guaranteed by \cite[Theorem 3.2]{FiaccoKyparisis85}.
Because the derivatives are finite a bound of the form $ \big \| \frac{\partial \b{a}}{ \partial \lambda_n } \big \|_2 \le K_a $
exists, which implies $ \big \| \frac{\partial \b{s}}{ \partial \lambda_n } \big \|_2 \le K_a $,
$ \big \| \frac{\partial \b{x}}{ \partial \lambda_n } \big \|_2 \le K_a $ \jbbr{and $ | \partial_n v_i | \le K_a $, $ | \partial_n \theta_i | \le K_a $.
Having deduced bounds on the sensitivities $ \partial_n v_i $ and $ \partial_n \theta_i $ now the partial derivatives of 
$ \frac{\partial \b{f}}{ \partial \b{x} } $ can be analyzed.}
\jbbr{Specifically}, the partial derivatives of $ \frac{\partial }{\partial \lambda_n} \left( \frac{\partial \b{f}}{ \partial \b{x} } \right) $ are
computed from \eqref{eq:jacelements} from which 
$ \frac{\partial}{\partial \lambda_n}\left( \frac{\partial \b{f}}{ \partial \biidxex{q}{G}{g} } \right) = \b{0} $. 
Since the power flow equations hold with equality at \jbbr{a} solution, we simplify the summations
in \eqref{eq:jacelements}, using the definitions $ p_i^{\text{net}} \equiv p_i^g - p_i^d $ $ q_i^{\text{net}} \equiv q_i^g - q_i^d $, as e.g., $ \sum_{k=1}^N v_k c_{ik} = \frac{p_i^{\text{net}}}{v_i} $
(and correspondingly for the other terms). Then for
$ 1 \le i \le N $ and $ j = N+ i $
as well as $ 1 \le l \le N_L $ and $ 1 \le t \le N $ 
\begin{align*}
	\frac{\partial}{\partial \lambda_n}\left( \frac{\partial \b{f}}{ \partial \biidx{v}{L} } \right)_{i,l} & = 
	\begin{cases}
	\partial_n (  \frac{p_i^{\text{net}}}{v_i} ) + 2\partial_n (c_{ii} v_i)	& \text{ if }  i = [L]_l  \\
		\partial_n (v_i c_{i[L]_l})  & \text{ otherwise }
	\end{cases}  \\
	 \frac{\partial}{\partial \lambda_n}\left( \frac{\partial \b{f}}{ \partial \biidx{v}{L} } \right)_{j,l} &=
	\begin{cases}
		\partial_n (  \frac{q_i^{\text{net}}}{v_i} ) + 2\partial_n (d_{ii} v_i) 	& \text{ if }  i = [L]_l  \\
		\partial_n (v_i d_{i[L]_l})					& \text{ otherwise }
	\end{cases}	\\	
	\frac{\partial}{\partial \lambda_n}\left( \frac{\partial \b{f}}{ \partial \bs{\theta} } \right)_{i,t} & = 
	\begin{cases}
		- \partial_n q_i^{\text{net}} 	 & \text{ if }  i = t  \\
		v_{k} d_{ik} \partial_n v_i+ v_i \partial_n (v_{k} d_{ik})				& \text{ otherwise }
	\end{cases}  \\
	\frac{\partial}{\partial \lambda_n} \left( \frac{\partial \b{f}}{ \partial \bs{\theta} } \right)_{j,t} &=
	\begin{cases}
		 \partial_n q_i^{\text{net}} 	& \text{ if }  i = t  \\
		-(v_{k} c_{ik} \partial_n v_i+ v_i \partial_n (v_{k} c_{ik})) 				& \text{ otherwise } \\
	\end{cases}  
\end{align*}
These derivatives can be bounded with limits on 
$ \partial_n p^{\text{net}}_i, \partial_n q^{\text{net}}_i, \partial_n v_i, \partial_n \theta_i, \partial_n c_{ik}, \partial_n d_{ik}, \partial_n (v_k c_{ik})  $, and $\partial_n (v_k d_{ik})$ \jbbo{(i.e., all elements that appear in the expression of the derivatives)}. 
\jbbo{Since all these individual elements are bounded by $ K_a $, and since the $ c_{ik}, d_{ik} $ terms are bounded, too 
(they depend on $ v_i, \theta_i $), the maximum element will be bounded by a constant, as well. Denote this constant
by $ K_x $, then
\begin{equation*}
	\underset{ 1\le i,j \le 2N }{\text{max}} \left | \left( \frac{\partial }{\partial \lambda_n} \left( \frac{\partial \b{f}}{ \partial \b{x} } \right)  \right)_{ij}  \right |
	\equiv K_x.
\end{equation*}}
With this, \jbbo{since $ \| \b{j}_r \|_2 \le K_{\Gamma} $ and 
$ \frac{\partial^2 \b{f}}{ \partial \lambda_n  \partial \b{x} } = \frac{\partial}{ \partial \lambda_n }  \b{J}(\b{s}(\bs{\lambda})) $}, we obtain the upper bound 
\begin{align}
	\label{eq:boundJac}
	\left\| \b{j}_r \tp \frac{\partial  \b{J}(\b{s}(\bs{\lambda})) }{ \partial \lambda_n }\right \|_2 \le
	2 K_{\Gamma} K_x N.
\end{align}
\subsection{Part III: Bound on $ \left\| \frac{\partial \bs{\lambda}_M}{ \partial \bs{\lambda} } \right \|_2 $  and Implications}
\label{sec:impl}
The analysis from the previous section has implications for the convergence of Algorithm 1. In 
particular, since $ \| \b{J}(\b{s}(\bs{\lambda}))^{-1} \|_2 \le K_{\Gamma} $ 
(cf. Section \ref{sec:derivsxo}) a bound on the magnitude of $ \frac{\partial \bs{\lambda}_M}{\partial \bs{\lambda}} $ can be found. By combining \eqref{eq:boundjacel} and \eqref{eq:boundJac} 
\jbbr{one sees that
\begin{align*}
			\left(\frac{\partial \bs{\lambda}_M}{\partial \bs{\lambda}} \right)_{rn} &\leq | z_r |~\left\| \bs{\Sigma} \right\|_2~\left\| \b{e}\tp_r (\partial \bs{\Gamma} \big /\partial \lambda_n)  \right\|_2 \\ 
			&\leq | z_r |~\left\| \bs{\Sigma} \right\|_2~K_{\Gamma}~\left\| \b{j}\tp_r (\partial \b{J} \big / \partial \lambda_n) \right\|_2 \\
			&\leq | z_r |\jbbR{\cdot} \left\| \bs{\Sigma} \right \|_2\jbbR{\cdot}K_{\Gamma}\jbbR{\cdot}2\jbbR{\cdot}K_{\Gamma}\jbbR{\cdot}K_x\jbbR{\cdot}N,
\end{align*}}
\jbbr{where the second inequality is based on \eqref{eq:gammaDerivBound} and the third on \eqref{eq:boundJac}.}
Since typically only $ N_A \ll 4N + N_L $ inequality constraints are active at the solution, \jbbr{and since $ | z_r | \le K_1 $} the upper bound
is 
\begin{equation}
	\label{eq:upperbound}
		\left \| \frac{\partial \bs{\lambda}_M}{\partial \bs{\lambda}} \right\|_2 \le 2 \cdot \| \bs{\Sigma} \|_2 \cdot K_1  \cdot K^{2}_{\Gamma} \cdot K_x \cdot N_A \cdot N.
\end{equation}
If the right hand side in \eqref{eq:upperbound} is less than one then Condition 2 in \eqref{eq:convcond} of the 
fixed point iteration is guaranteed to hold. \jbbR{In other words
\begin{align}
	\label{eq:uppReform}
	0 \le 2 \cdot \| \bs{\Sigma} \|_2 \cdot K_1 & \cdot K^{2}_{\Gamma} \cdot K_x \cdot N_A \cdot N <1 \quad \text{implies} \nonumber  \\ 
	&\text{Algorithm 1 converges}. 
\end{align}}
\subsection{Computing a Bound Estimate}
\label{sec:practice}
The upper bound in \eqref{eq:upperbound} is composed of the
parts
\begin{equation*}
	\| \bs{\Sigma} \|_2, \quad K_1, \quad  K^2_{\Gamma}  \cdot  K_x \cdot N_A \cdot N.
\end{equation*}
\jbbR{The quantities $\|\bs{\Sigma}\|_2$ and $K_1$ are related to the uncertainty in the model,}
whereas $ K^2_{\Gamma} \cdot K_x \cdot N_A \cdot N $ 
is problem (network specification) dependent. For instance, $K_1$ is reduced by reducing the probability constraints
in \eqref{eq:ccacopf}. Concretely, if all $ \epsilon \to 1/2 $ then $K_1 \to 0$. Secondly, \jbb{\emph{if $ \| \b{\Sigma} \|_2 $ is small enough (sufficiently small uncertainty) then  from \eqref{eq:convcond} the fixed point iteration is guaranteed to converge}}.
\jbbr{Although the bounds may be known to be conservative,
	we show that if the uncertainty is small enough convergence is achieved. Such an insight may be useful if the
	uncertainty can be scaled (e.g., using shorter model horizons) so that its magnitude becomes smaller. 
	Such conclusions apply to basically any type of Algorithm 1 as long as $ \boldsymbol{\Sigma} $
	is used in forming $ \boldsymbol{\lambda} $.}
Finally, the bound \eqref{eq:upperbound} includes the network dimension $N$. 
In numerical experiments, we observe that setting $ \| \bs{\Sigma} \|_2 \sim \mathcal{O}(1/N^2) $ enables the 
solution of large networks.
Before applying Algorithm 1 for multiple iterations the bound in \eqref{eq:upperbound} can be qualitatively
used to asses the convergence behavior of the method. In particular, at $k=0$, use $ \b{s}^{(1)} $
to compute \jbbr{parts of} the right hand side in \eqref{eq:upperbound}. Computing $ K_{\Gamma} $ 
is available from \jbbr{or $ K_{\Gamma} = \sqrt{\| \bs{\Gamma} \|_1 \| \bs{\Gamma} \|_{\infty}}  $ (or a LU factorization of $ \b{J} $)}. 
\jbbo{The computation of $ K_x $ can be approximated
by first estimating $ K_a $ from a lower bound on the smallest singular value 
of  $ \partial \b{G} / \partial \b{a} $ using \cite{HongPan92} (as for \eqref{eq:jacbound}), calling such an estimate $ \hat{K}_a $. Because $ \bs{\lambda} $
appears linearly in the inequalities $ \b{h}(\b{s},\bs{\lambda}) =\text{vec}(\b{g}_1(\b{s}), \b{g}_2(\b{s}) + \bs{\lambda})$ (for appropriately defined $\b{g}_1(\b{s}), \b{g}_2(\b{s}) $)
and in no other constraints,
the expression $ \partial \b{G} / \partial \lambda_n = \jbb{\text{vec}(\b{0},\b{e}_n)} $ simplifies. Subsequently, $ K_a \sim 1 / \hat{K}_a $
and $ K_x = a K_a $, for a constant $a > 0$. \jbbr{Because the estimate for $ K_x $ may be reasonably small near a solution, and because $ K_x $ would incur extra computational expenses we set its value to $ K_x \in (0,1] $. } } 
If the computed bound, say $ B^{(0)} $ exceeds a 
fixed threshold (say $ \text{threshold} = 10 $), set $ \bs{\Sigma}^{\text{Scaled}} = 1 / B^{(0)} \bs{\Sigma} $. 
Moreover, we include a scaling factor $ \gamma_g $ in computing the line flow constraint
tightenings from \eqref{eq:lamconsLF}, thus $ \bar{\lambda}_{r_1} = \gamma_g \lambda_{r_1} $.
Note that for $ \gamma_g = 1 $ no scaling is included.
\jbbr{An overview of the constants is summarized in a table. Furthermore, the quantity $ K_P = \| \bs{\Sigma} \| K^2_{\Gamma} N_A $ captures
specific problem sensitivities to relevant factors such as magnitude of $ \bs{\Gamma} $ and number of active constraints. This quantity can be computed
at the beginning of applying Algorithm 1, is often inexpensive, and is expected to be small for guaranteed convergence.}
\begin{table}[h!]
\label{tbl:constants}
\caption{\jbbr{Components in computing a bound estimate in \eqref{eq:upperbound}}}
\centering
\begin{tabular}{ | c | p{3.2cm} | p{3.2cm} | }
\hline \textbf{Constant} 			& \textnormal{Meaning} 						& \textnormal{Computation} \\
\hline$ K_1 $ 					& Distribution bound ($ z_r = F^{-1}(1-\epsilon_r) $)		& $\underset{r}{\textnormal{max}} | z_r | $ \\
$ K_{\Gamma} $ 				& Magnitude inv. power flow Jac. 				&  $ \sqrt{\| \bs{\Gamma} \|_1 \| \bs{\Gamma} \|_{\infty}} $ or eq. \eqref{eq:jacbound} \\
$ K_{x} $ 						& Magnitude of power flow Jac. deriv. w.r.t $ \lambda $ 	& $\in (0,1]$ \\
$ N_A $ 						& Number of active inequality constraints (w/o deter. vars)			& $\text{count}(\b{h}(\b{s}^{(1)})=0) $ \newline (w/o constraints on $ \b{y} $) \\
$ K_P $						& Problem sensitivity to inv. Jac. and active constraints 	& $\| \bs{\Sigma} \| \cdot K^2_{\Gamma} \cdot N_A$ \\
 \hline  
\end{tabular}
\end{table}
\jbbr{(If $ \|\bs{\Sigma} \| $ is not known it can be bounded by $ K_2 = \sqrt{\| \bs{\Sigma} \|_1 \| \bs{\Sigma} \|_{\infty}} $)}

\section{Numerical Experiments}
\label{sec:numerical}
This section describes numerical experiments on a set of standard IEEE test problems. Typically,
solving large chance-constrained optimization problems is numerically very challenging. Even
directly solving the reformulated problem in \eqref{eq:iacopf} with state-of-the-art general
purpose methods becomes quickly intractable. Therefore, for larger networks we
apply Algorithm 1 to approximately solve \eqref{eq:iacopf}. Our implementation
of Algorithm 1 uses Julia 1.1.0 \cite{bezanson2017julia} and the modeling libraries
JuMP v0.18.6 \cite{dunning2017jump} and StructJuMP \cite{structJump}. In particular, in order to solve
the AC-OPF problem from \eqref{eq:acopf} and its modifications in Line \ref{alg:alg1_s} of Algorithm 1
we use the general purpose nonlinear optimization solver IPOPT \cite{WaechterBiegler06}. 
The stopping criteria in Algorithm 1 are set as $ \tau_q = 1 \times 10^{-3} $, $ \tau_v = 1 \times 10^{-5} $,
$ \tau_\theta = 1 \times 10^{-5} $, $ \tau_g = 1 \times 10^{-3} $. Unless otherwise stated, we set $ \bs{\Sigma} = \sigma \b{I} $, $ \sigma = 1/ N^2 $
and $ \epsilon = \text{vec}(\epsilon_q, \epsilon_v, \epsilon_\theta, \epsilon_g) = \text{vec}(0.1,0.1,0.1,0.2)  $, and $  \gamma_g = 1 / N_L^2 $.
The 
initial values are computed as the midpoint between the upper and lower bounds $ \b{s}^{(0)} = 1/2 (\b{u} + \b{l}) $. The maximum 
number of iterations for Algorithm 1 is set as \jbbR{50}. 
The test cases are summarized in Table II: 
\begin{table}[h!]
\label{tbl:cases}
\caption{Test Cases. Column 2 contains the total number 
of nodes in the network, $N$, and their split among generators, $N_G$, and loads, $N_L$. The largest
case contains more \jbbR{than $9000$ nodes}.} 
\centering
\begin{tabular}{ | l | l | l | c| }
\hline \textbf{Problem} 			& $ N = N_G + N_L $ 						& $ N_{\text{line}} $ & \text{Reference} \\
\hline \text{IEEE 9} 				& $9 = 3 + 6 $ 		 	 					& $ 9$ & \cite{matpower11} \\
\text{IEEE 30} 					& $30 = 6 + 24 $ 				 			& $41 $& \cite{matpower11,alsacstott74} \\
\text{IEEE 118} 					& $118 = 54 + 64 $ 				 		& $186 $&  \cite{matpower11} \\
\text{IEEE 300} 				& $300 = 69 + 231 $ 				  			& $411 $& \cite{matpower11} \\
\text{IEEE 1354pegase} 			& $1354 = 260 + 1094 $  						& $1991$ & \cite{FPCW13} \\
\jbbr{\text{IEEE 2383wp}} 			& $\jbbr{2383 = 327 + 2056} $  						& $\jbbr{2896}$ & \cite{matpower11} \\
\text{IEEE 2869pegase} 					& $2869 = 510 + 2359 $ 				& $4582$ &  \cite{FPCW13} \\
\text{IEEE 9241pegase} 					& $9241 = 1445 + 7796$ 				& $16049$ & \cite{FPCW13}  \\
 \hline  
\end{tabular}
\end{table}
The numerical experiments are divided into \jbbr{three} parts. Experiment I compares solutions to
\eqref{eq:iacopf} using a direct solver (i.e., IPOPT) or the iterative approach from Algorithm 1.
Experiment II reports results on the test problems from Table II.
\jbbr{Experiment III describes convergence tests according to our analysis, outcomes for problems with perturbed loads
and an investigation of joint and \jbbR{disjoint} relative frequencies on the IEEE 9 problem.} 

\subsection{Experiment I}
\label{subsec:EXI}
\begin{figure}[!t]
\centering
\includegraphics[width=2.5in]{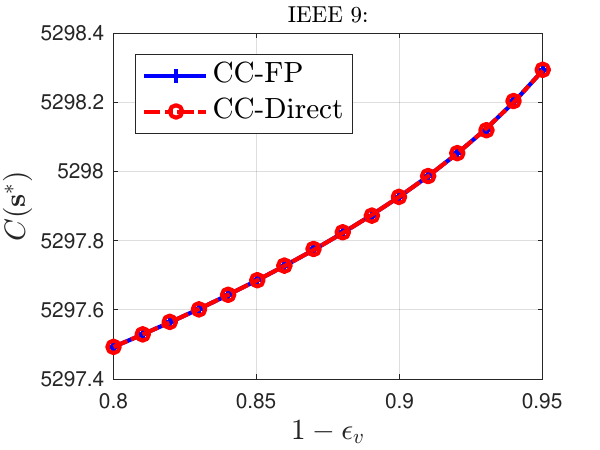}
\caption{Comparison of optimal objective function values using \jbbr{CC-Direct} (direct approach) and 
\jbbr{CC-FP (fixed-point)}. The optimal objective function values nearly coincide on this IEEE 9 bus case
as the parameters $ \bs{\epsilon}_v $  (probability thresholds) vary. }
\label{fig:EXIa}
\end{figure}

\begin{figure}[!t]
\centering
\includegraphics[width=2.5in]{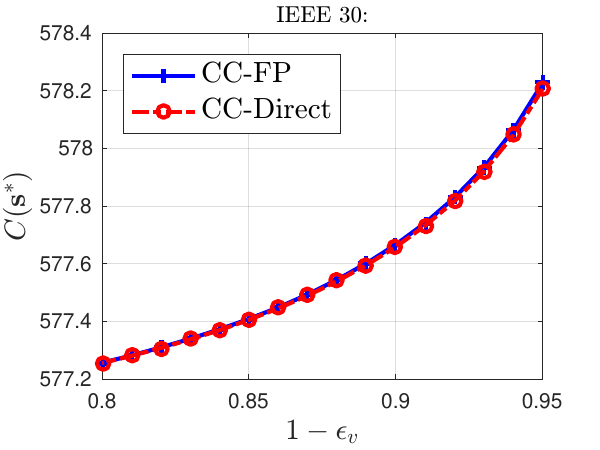}
\caption{Comparison of optimal objective function values using \jbbr{CC-Direct} (direct approach) and 
\jbbr{CC-FP (fixed-point)}. Also on this IEEE 30 bus case the optimal objective function values nearly coincide
as the parameters $ \bs{\epsilon}_v $  (probability thresholds) vary. }
\label{fig:EXI}
\end{figure}

This experiment compares the optimal objective function values of solving \eqref{eq:iacopf} directly
(\jbbr{CC-Direct}) using IPOPT or using Algorithm 1 \jbbr{with a fixed point approach} (\jbbr{CC-FP}). Solving
\eqref{eq:iacopf} in JuMP  becomes computationally intractable beyond $ N=100 $, because
the inverse, $ \bs{\Gamma}(\b{s}) \in \mathbb{R}^{2N \times 2N} $, and its derivatives are
continuously recomputed. For this reason, we compare \jbbr{CC-Direct} and \jbbr{CC-FP}
for the 9 bus and 30 bus cases (which are both tractable by the direct approach). 
\jbbr{For the purpose of this comparison we switched the line flow tightening off (i.e., $ \bs{\lambda}_g = \b{0} $) while all other
tightenings are applied}. In particular,
we compare the solved objective function values after perturbing the problem formulation
slightly. Specifically, we compare the computed objective function values when the values 
$ \epsilon_v = \{0.05,0.06,\ldots,0.2\} $ change the probability constraints.
(Varying other parameters yields similar results). The outcomes of the 9 bus network are in Figure \ref{fig:EXIa}.
Figure \ref{fig:EXI} contains the outcomes of applying both approaches on the 30 bus network.  
Observe in both figures that the optimal objective function values, i.e., $ C(\b{s}^*) $,
increase as the values of $ 1 - \epsilon_v $ increase. This is because larger values
of  $ 1 - \epsilon_v $ correspond to stricter chance constraints. Importantly, observe
that the objective function values are virtually equivalent for these test cases,
as the values of the blue and red curves nearly exactly match up. However, the main advantage of Algorithm 
1 is that it scales to larger cases, too.

\subsection{Experiment II}
\label{subsec:EXII}
Experiment II reports results of applying Algorithm 1 (\jbbr{CC-FP})
and a direct solver (\jbbr{CC-Direct}) to the problem from \eqref{eq:iacopf} on the test problems
in Table \jbbr{II}. For \jbbr{CC-FP}
the number of iterations, time and ``optimal" objective values are reported. For \jbbr{CC-Direct} the
``optimal" objective values and the time are reported. Algorithm 1 converged 
on all of the reported problems. 
These problems include large network instances, too.
Observe in Table \jbbr{III} that for problems on which both solvers can be used, the optimal
objective function values are close to each other. However, for large problems only Algorithm 1,
based on a fixed point iteration, is practical. \jbbr{For consistency with Figs. \ref{fig:EXIa} and \ref{fig:EXI} the line flow tightenings are switched off
in the first two cases, while all tightenings are applied on all other problems.}

\begin{table}
	\label{tbl:EXII}
	\caption{ Comparison of Algorithm 1 (\jbbr{CC-FP}) and 
	a direct solution approach (\jbbr{CC-Direct}) on a set of IEEE power networks. The last 6 rows
	in the ``Objective" column for \jbbr{CC-Direct} are labeled $ \texttt{N/A} $, because this approach
	took computational times in excess of $5$ hours.}
	\small
	\begin{tabular}{ | p{0.7cm} | p{0.25cm} | c | c | c | p{1.05cm} | }
		\hline
		\multirow{3}{*}{IEEE}   
		& \multicolumn{3}{c|}{\jbbr{CC-FP}} 
		& \multicolumn{2}{c|}{\jbbr{CC-Direct}}\\
		\cline{2-6}
		&  Its.
		& \jbbr{Obj. (Cost/h)}
		& Time(\jbbr{s})
		& \jbbr{Obj. (Cost/h)}
		& Time(\jbbr{s}) \\ 
		\hline
		$9$  & $\texttt{4}$&   				$\texttt{5297.928}$&   $\texttt{0.0504}$&   $\texttt{5297.928}$&   $\texttt{0.3075}$\\
		$30$ & $\texttt{4}$&   				$\texttt{577.6665}$&   $\texttt{0.1688}$&   $\texttt{577.6592}$&   $\texttt{148.36}$\\ 
		$118$ & $\texttt{3}$&   			$\texttt{129662.0}$&   $\texttt{0.4710}$&   $\texttt{N/A}^{\dagger}$&   $\texttt{>5h}$\\
		$300$ & $\texttt{5}$&   			$\texttt{720090.3}$&  $\texttt{4.3123}$&    $\texttt{N/A}^{\dagger}$&   $\texttt{>5h}$\\
		$1354$ & $\texttt{3}$&   	$\texttt{74069.38}$&   $\texttt{75.039}$&   $\texttt{N/A}^{\dagger}$&   $\texttt{>5h}$\\
		$2383$ & $\texttt{3}$&   	$\texttt{1868551.}$&   $\texttt{236.27}$&   $\texttt{N/A}^{\dagger}$&   $\texttt{>5h}$\\
		$2869$ & $\texttt{3}$&   	$\texttt{133999.3}$&   $\texttt{316.80}$&   $\texttt{N/A}^{\dagger}$&   $\texttt{>5h}$\\
		$9241$ & $\texttt{3}$&   	$\texttt{315912.6}$&   $\texttt{3255.2}$&   $\texttt{N/A}^{\dagger}$&   $\texttt{>5h}$\\
		\hline
	\end{tabular}
\end{table}


\subsection{\jbbr{Experiment III}}
\label{subsec:EXIII}
\jbbr{In this experiment we investigate three further important aspects of the analysis and the algorithm. 
\jbbR{Before proceeding, we note that in Algorithm 1 it can happen that the constraint tightenings become
large enough that adjusted upper bound constraints fall below a lower bound constraint (or vice-versa). Since this situation
implies an inconsistent constraint, we restore consistency by re-setting this particular constraint to its original,
yet with a scaled feasible interval (to reflect some tightening). Use of this correction mechanism increases the robustness
of the method, meaning that Algorithm 1 converges more often. Nevertheless, our analysis indicates that
a sufficiently small value of $ \| \bs{\Sigma} \| $ will guarantee convergence. However, other network dependent factors
such as $ \| \bs{\Gamma} \| $ or $ N_A $ (number of active constraints) are also relevant for convergence on specific
problems. We capture the conclusion that a sufficiently small $ \| \bs{\Sigma} \| $ guarantees convergence, and that
relatively large values of $ \|  \bs{\Gamma} \| $ and $N_A$ typically result in non-convergence using all
test problems from Table II.}} 

\subsubsection{\jbbr{Experiment III.A}}
\label{subsubsec:EXIIIa}
\jbbr{This experiment varies the default value $ \sigma = 1 / N^2 $ as $ \sigma = \alpha / N^2 $ with
$ \alpha = \{1, 10^2, 10^4, 10^6, \jbbR{10^8} \} $. For a sufficiently small value of $ \sigma $ (here $\alpha=1$) all
problems are guaranteed to converge. Moreover, lower values of $ \| \bs{\Gamma} \| $ and $N_A$ represent
lower sensitivity of a problem to increases in $ \| \bs{\Sigma} \| $. Table~IV reports the 
outcomes on all test problems.} 
\jbbR{This table contains explicit values for all parameters from Table III. Thus the bound estimate from \eqref{eq:upperbound} can
be computed. Particularly, note that when $\alpha=1$ all the bounds are $<1$ (cf. the row labeled \texttt{Cond. 30}). Moreover,
when $\alpha=10$ all problems with bound estimate $<1$ converged (i.e., {\small IEEE \{9,118,1354,2869,9241\}}). This is consistent with \eqref{eq:uppReform}.
Since Condition 2 in \eqref{eq:convcond} is sufficient and our analysis develops an upper bound, networks with bound estimate larger than 1 many also converge.
However, when the upper bound is less than 1 then convergence of Algorithm 1 is implied.
	Additionally, problems {\small IEEE $\{1354, 2869\}$} have the smallest 
	$ K_P $ values indicating their tolerance to larger $\sigma$ values. Problem {\small IEEE $300$} has the largest 
	$ K_P $ value and is viewed as being more sensitive to non-convergence. Ultimately 
	all problems are non-convergent. Prior to this {\small IEEE $\{1354, 2869\}$} tend to converge on more instances while these networks have the smallest 
	$ K_P $ values.}

\begin{table*}
\centering
\label{tbl:EXIII_UP}
	\caption{\jbbR{Convergence of Algorithm 1 when $ \sigma = \frac{\alpha}{N^2} $
	for $ \alpha = \{1, 10, 10^4, 10^6, 10^8 \} $. A boxed $\texttt{Y}$ in rows labeled~\texttt{Conv.}~indicates convergence. Values in rows labeled 
	~\texttt{Cond. 30}~correspond to the computed value in the right hand side of \eqref{eq:upperbound}. Bold or 
	resp. italic $ K_P $ values represents the smallest resp. second smallest $K_P$ values.}}
\begin{tabular}{| c | c | c c c | l | c c c c c |}
	\hline
\multicolumn{1}{| c |}{} & \multirow{2}{*}{IEEE} & \multicolumn{3}{c |}{} & & \multicolumn{5}{c |}{ $ \alpha $}\\  \cline{7-11} 
$K_1$ & & $K_x$  & $K_{\Gamma}$ & $N_A$ & & $1$ & $10$ & $10^4$ & $10^6$ & $10^8 $ \\ 
\hline
\multirow{32}{*}{$1.3$} &\multirow{4}{*}{$9$} &\multirow{4}{*}{$0.005$} &\multirow{4}{*}{$1.6$} &\multirow{4}{*}{$2$} & $\sigma$ & $\texttt{0.012}$&$\texttt{0.12}$&$\texttt{1.2e+02}$&$\texttt{1.2e+04}$&$\texttt{1.2e+06}$\\ 
& & & & &  $K_{P}$ & $\texttt{0.063}$&$\texttt{0.63}$&$\texttt{6.3e+02}$&$\texttt{6.3e+04}$&$\texttt{6.3e+06}$\\ 
& & & & &  $\texttt{Conv.}$ & ${\boxed{\texttt{Y}} }$&${\boxed{\texttt{Y}} }$&${\boxed{\texttt{Y}} }$&${ \texttt{N}}$&${ \texttt{N}}$\\ 
& & & & &  $\texttt{Cond. 30}$ & $\texttt{0.0072}$&$\texttt{0.072}$&$\texttt{72}$&$\texttt{7.2e+03}$&$\texttt{7.2e+05}$\\\cline{2-11} 
&\multirow{4}{*}{$30$} &\multirow{4}{*}{$0.1$} &\multirow{4}{*}{$8.5$} &\multirow{4}{*}{$1$} & $\sigma$ & $\texttt{0.0011}$&$\texttt{0.011}$&$\texttt{11}$&$\texttt{1.1e+03}$&$\texttt{1.1e+05}$\\ 
& & & & &  $K_{P}$ & $\texttt{0.08}$&$\texttt{0.8}$&$\texttt{8e+02}$&$\texttt{8e+04}$&$\texttt{8e+06}$\\ 
& & & & &  $\texttt{Conv.}$ & ${\boxed{\texttt{Y}} }$&${ \texttt{N}}$&${ \texttt{N}}$&${ \texttt{N}}$&${ \texttt{N}}$\\ 
& & & & &  $\texttt{Cond. 30}$ & $\texttt{0.62}$&$\texttt{6.2}$&$\texttt{6.2e+03}$&$\texttt{6.2e+05}$&$\texttt{6.2e+07}$\\\cline{2-11} 
&\multirow{4}{*}{$118$} &\multirow{4}{*}{$0.005$} &\multirow{4}{*}{$6.3$} &\multirow{4}{*}{$16$} & $\sigma$ & $\texttt{7.2e-05}$&$\texttt{0.00072}$&$\texttt{0.72}$&$\texttt{72}$&$\texttt{7.2e+03}$\\ 
& & & & &  $K_{P}$ & $\texttt{0.045}$&$\texttt{0.45}$&$\texttt{4.5e+02}$&$\texttt{4.5e+04}$&$\texttt{4.5e+06}$\\ 
& & & & &  $\texttt{Conv.}$ & ${\boxed{\texttt{Y}} }$&${\boxed{\texttt{Y}} }$&${ \texttt{N}}$&${ \texttt{N}}$&${ \texttt{N}}$\\ 
& & & & &  $\texttt{Cond. 30}$ & $\texttt{0.069}$&$\texttt{0.69}$&$\texttt{6.9e+02}$&$\texttt{6.9e+04}$&$\texttt{6.9e+06}$\\\cline{2-11} 
&\multirow{4}{*}{$300$} &\multirow{4}{*}{$0.005$} &\multirow{4}{*}{$23$} &\multirow{4}{*}{$31$} & $\sigma$ & $\texttt{1.1e-05}$&$\texttt{0.00011}$&$\texttt{0.11}$&$\texttt{11}$&$\texttt{1.1e+03}$\\ 
& & & & &  $K_{P}$ & $\texttt{0.19}$&$\texttt{1.9}$&$\texttt{1.9e+03}$&$\texttt{1.9e+05}$&$\texttt{1.9e+07}$\\ 
& & & & &  $\texttt{Conv.}$ & ${\boxed{\texttt{Y}} }$&${ \texttt{N}}$&${ \texttt{N}}$&${ \texttt{N}}$&${ \texttt{N}}$\\ 
& & & & &  $\texttt{Cond. 30}$ & $\texttt{0.73}$&$\texttt{7.3}$&$\texttt{7.3e+03}$&$\texttt{7.3e+05}$&$\texttt{7.3e+07}$\\\cline{2-11} 
&\multirow{4}{*}{$1354$} &\multirow{4}{*}{$0.02$} &\multirow{4}{*}{$5.1$} &\multirow{4}{*}{$90$} & $\sigma$ & $\texttt{5.5e-07}$&$\texttt{5.5e-06}$&$\texttt{0.0055}$&$\texttt{0.55}$&$\texttt{55}$\\ 
& & & & &  $K_{P}$ & $\texttt{\emph{0.0013}}$&$\texttt{\emph{0.013}}$&$\texttt{\emph{13}}$&$\texttt{\emph{1.3e+03}}$&$\texttt{\emph{1.3e+05}}$\\ 
& & & & &  $\texttt{Conv.}$ & ${\boxed{\texttt{Y}} }$&${\boxed{\texttt{Y}} }$&${\boxed{\texttt{Y}} }$&${\boxed{\texttt{Y}} }$&${ \texttt{N}}$\\ 
& & & & &  $\texttt{Cond. 30}$ & $\texttt{0.089}$&$\texttt{0.89}$&$\texttt{8.9e+02}$&$\texttt{8.9e+04}$&$\texttt{8.9e+06}$\\\cline{2-11} 
&\multirow{4}{*}{$2383$} &\multirow{4}{*}{$0.001$} &\multirow{4}{*}{$45$} &\multirow{4}{*}{$253$} & $\sigma$ & $\texttt{1.8e-07}$&$\texttt{1.8e-06}$&$\texttt{0.0018}$&$\texttt{0.18}$&$\texttt{18}$\\ 
& & & & &  $K_{P}$ & $\texttt{0.092}$&$\texttt{0.92}$&$\texttt{9.2e+02}$&$\texttt{9.2e+04}$&$\texttt{9.2e+06}$\\ 
& & & & &  $\texttt{Conv.}$ & ${\boxed{\texttt{Y}} }$&${\boxed{\texttt{Y}} }$&${\boxed{\texttt{Y}} }$&${ \texttt{N}}$&${ \texttt{N}}$\\ 
& & & & &  $\texttt{Cond. 30}$ & $\texttt{0.56}$&$\texttt{5.6}$&$\texttt{5.6e+03}$&$\texttt{5.6e+05}$&$\texttt{5.6e+07}$\\\cline{2-11} 
&\multirow{4}{*}{$2869$} &\multirow{4}{*}{$0.01$} &\multirow{4}{*}{$5.9$} &\multirow{4}{*}{$185$} & $\sigma$ & $\texttt{1.2e-07}$&$\texttt{1.2e-06}$&$\texttt{0.0012}$&$\texttt{0.12}$&$\texttt{12}$\\ 
& & & & &  $K_{P}$ & $\texttt{\textbf{0.00078}}$&$\texttt{\textbf{0.0078}}$&$\texttt{\textbf{7.8}}$&$\texttt{\textbf{7.8e+02}}$&$\texttt{\textbf{7.8e+04}}$\\ 
& & & & &  $\texttt{Conv.}$ & ${\boxed{\texttt{Y}} }$&${\boxed{\texttt{Y}} }$&${\boxed{\texttt{Y}} }$&${\boxed{\texttt{Y}} }$&${ \texttt{N}}$\\ 
& & & & &  $\texttt{Cond. 30}$ & $\texttt{0.057}$&$\texttt{0.57}$&$\texttt{5.7e+02}$&$\texttt{5.7e+04}$&$\texttt{5.7e+06}$\\\cline{2-11} 
&\multirow{4}{*}{$9241$} &\multirow{4}{*}{$0.002$} &\multirow{4}{*}{$18$} &\multirow{4}{*}{$511$} & $\sigma$ & $\texttt{1.2e-08}$&$\texttt{1.2e-07}$&$\texttt{0.00012}$&$\texttt{0.012}$&$\texttt{1.2}$\\ 
& & & & &  $K_{P}$ & $\texttt{0.002}$&$\texttt{0.02}$&$\texttt{20}$&$\texttt{2e+03}$&$\texttt{2e+05}$\\ 
& & & & &  $\texttt{Conv.}$ & ${\boxed{\texttt{Y}} }$&${\boxed{\texttt{Y}} }$&${\boxed{\texttt{Y}} }$&${ \texttt{N}}$&${ \texttt{N}}$\\ 
& & & & &  $\texttt{Cond. 30}$ & $\texttt{0.093}$&$\texttt{0.93}$&$\texttt{9.3e+02}$&$\texttt{9.3e+04}$&$\texttt{9.3e+06}$\\\cline{1-11} 
\end{tabular}
\end{table*}


\subsubsection{\jbbr{Experiment III.B}}
\label{subsub:EXIIIb}

\begin{figure}[!t]
\centering
\includegraphics[width=2.5in]{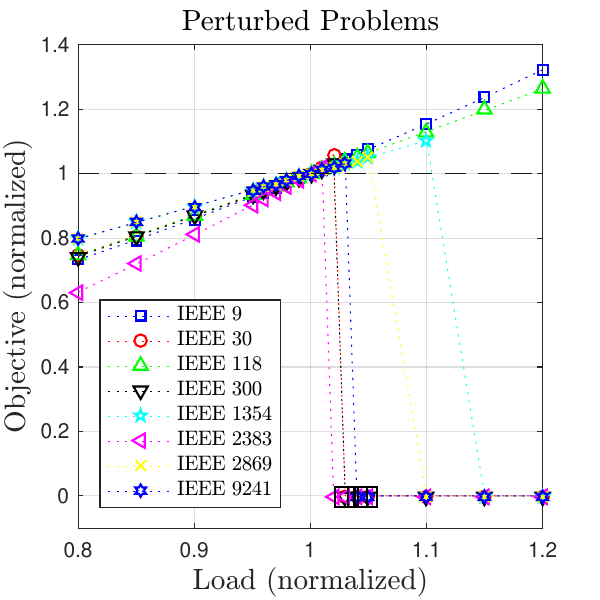}
\caption{\jbbr{Normalized optimal objective values (relative to the original optimal value), when nonzero 
loads at all buses in all systems vary by a factor of $0.8$ --- $1.2$. Objective values of $0$ indicate that the 
problem did not converge. All problems were solvable with load values up to one percent of the original. Some problems
such as {\small IEEE 9 } and {\small IEEE \jbbR{118}} were solvable with increases of up to $ 20\% $. Black squares represent
problems in which the corresponding classical ACOPF problem (without a probability model) did converge on the
same perturbed problem. In only 4 instances out of 136 total perturbed problems did this occur. Overall,
the changes in loads linearly affect the optimal values and the chance-constrained models have similar
sensitivities to load changes as the classical ACOPF.}}
\label{fig:PerturbedProbs}
\end{figure}

\jbbr{This experiment investigates outcomes on all test problems 
when loads are perturbed (representing stressed system conditions). In particular, the algorithm is applied
to all problems with loads varying between $ 80\% $ --- $120\% $ of their original values. Figure \ref{fig:PerturbedProbs} 
displays the outcomes in which optimal objective values relative to the original optimal value are compared.
Overall, the optimal objective values exhibit a linear relation with changes in the loads.}

\subsubsection{\jbbr{Experiment III.C}}
\label{subsub:EXIIIc}

\begin{figure}[!t]
\centering
\includegraphics[width=2.5in]{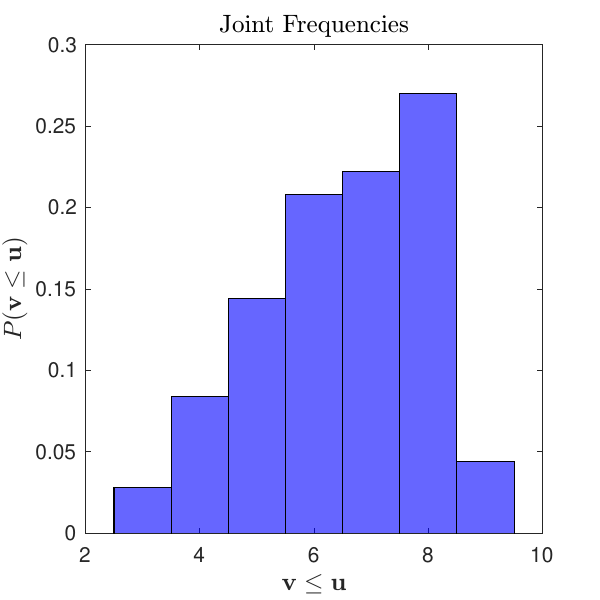} \\
\includegraphics[width=2.5in]{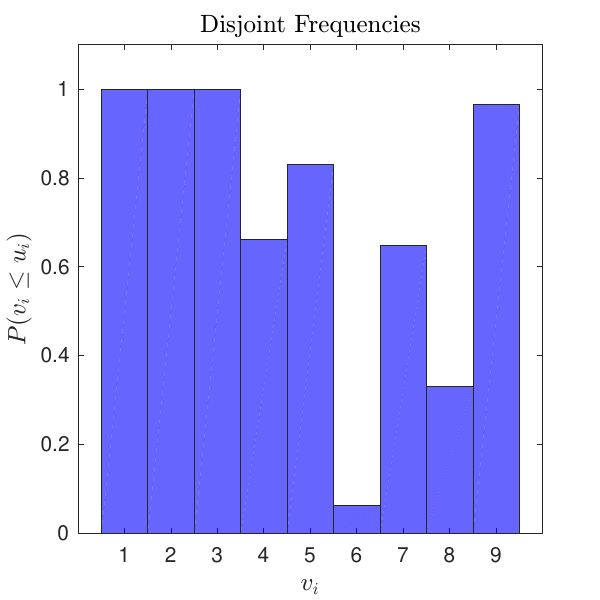} \\
\caption{\jbbr{Top: Relative frequencies for the count of elements of vector $ \b{v} $
				that satisfy $ \b{v} \le \b{u} $. The constraint that all 9 voltages simultaneously satisfy the inequality
				has probability of 0.044. 
			Bottom: Relative frequencies of each individual voltage to satisfy its constraint. Independently enforcing
			each voltage to satisfy the constraints $ \prod_{i=1}^9 P(v_i \le u_i) $ has a probability of 0.007.
			In this case there is a factor of about 6 (i.e. $ 0.044 / 0.007 $) between enforcing a joint or individual
			constraints. However, enforcing joint probability constraints as part of an optimization for large problems
			is typically not computationally feasible.}}
\label{fig:frequencies}
\end{figure}

\jbbr{This experiment investigates the relation of joint and individual (disjoint) probabilities. For the
{\small IEEE 9} bus case, we sample 500 multivariate random vectors $ \bs{\omega}^p, \bs{\omega}^q $ with
zero mean and positive definite dense covariance. The voltage constraints on generator
buses (for this 9 bus case i.e., 1,2 and 3) are deterministic $ v_j \le u_j, j=1,2,3 $. The voltages at the remaining
load buses are unbounded. The problem is re-solved after random loads $ \b{p}^d + \bs{\omega}^p $
and $ \b{q}^d + \bs{\omega}^q $ are simulated. The joint relative frequencies of the vector constraint 
$ \b{v} \le \b{u} $ (the constant vector $\b{u} $ has elements $ u_i = 1.1 $) are displayed on the top of
Figure \ref{fig:frequencies}. The individual relative frequencies for each bus are at the bottom of the figure.}

\section{Conclusion}
\jbbR{We briefly summarize limitations and future work of the analysis before concluding.}
\subsection{\jbbR{Limitations}}
\jbbR{Recall that Condition 2 in \eqref{eq:convcond} is a sufficient condition for convergence of a fixed point iteration
(such as Algorithm 1). This means that even when the condition does not hold the algorithm may still converge. In particular,
when the upper bound in \eqref{eq:upperbound} exceeds one the algorithm may still converge. However, when \eqref{eq:upperbound}
is less than one and assumptions A.1--A.3 from Sec. III.A are satisfied then convergence is guaranteed by the fixed point theorem. Moreover,
the components of \eqref{eq:upperbound} are computed using upper bounds or estimates as described in Table II and may thus be conservative.} 
\subsection{\jbbR{Future Work}}
\jbbR{Since the analysis is general, it applies basically to all fixed point algorithms which 
iteratively re-solve power optimization problems parameterized by a sequence of constraint tightenings, $\boldsymbol{\lambda}^{(k)}$.
Therefore, in future work one can apply the analysis to the convergence of fixed point approaches for variations such as Robust AC Optimal Power flows, High Voltage Direct Current lines or similar problems. The 
decomposition of the upper bound in \eqref{eq:upperbound} is already detailed. Yet, new work can tighten this upper bound with sharped estimates
for the constants, and particularly the constant $K_x$.}

In sum, this article describes a chance constrained AC optimal 
power flow model with only deterministic variables in the objective, which thus enables immediate
interpretation of the optimal function values. By linearizing
the stochastic variables, a deterministic nonlinear optimization problem
is obtained in lieu of the probabilistic one. Because solving
the reformulated optimization problem is computationally 
challenging, we \jbbR{develop} and  analyze the convergence
criteria for a fixed point algorithm that solves a sequence of 
modified AC optimal power flow problems and scales
to larger network sizes. The analysis connects the variance
of uncertainty and constraint probabilities to the 
convergence properties of the algorithm. In numerical experiments,
we compare the fixed point iteration to directly solving
the reformulated problem and test the method on IEEE
problems, including a network with over 9000 nodes.
\jbb{Certainly our bounds are quite conservative in this version, however this is the first 
attempt at proving convergence of the approach.}
\jbbR{This opens up future work, since iteratively resolving a sequence of tractable optimization problems
(by holding specific nonlinear terms fixed) has been very effective on a variety of power systems
applications.} 


\section*{\jbb{Acknowledgments}}
This work was supported by the U.S. Department of Energy, Office of Science,
Advanced Scientific Computing Research, under Contract DE-AC02-06CH11357
at Argonne National Laboratory and by NSF through award CNS-51545046.
We acknowledge fruitful discussion with Dr. Line Roald, who pointed out
that Algorithm 1 was observed to cycle between different points, and encouraged
an analysis of its convergence criteria. We also thank Eric Grimm, who helped in carrying
out parts of the numerical experiments in Section IV. \jbbr{Careful reading and helpful suggestions
of the editor and three anonymous reviewers markedly improved the manuscript.}

\ifCLASSOPTIONcaptionsoff
  \newpage
\fi



\bibliographystyle{IEEEtran}
\bibliography{myrefs}
%

%
%
%
%

%

\begin{IEEEbiography}[{\includegraphics[width=1.0in,height=1.25in,clip,keepaspectratio]{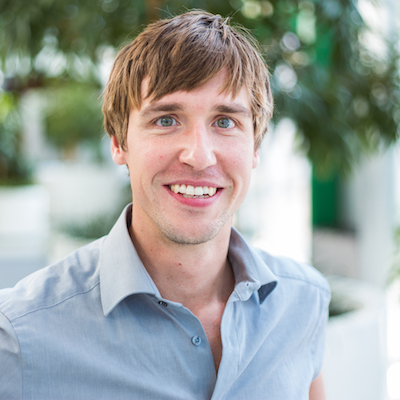} \vspace{0.2cm}}]{Johannes J. Brust} 
Was with the Mathematics and Computer Science Division at Argonne National
Laboratory, IL (now Department of Mathematics, University of California, San Diego). He received a M.Sc. in financial engineering from Maastricht University and a Ph.D. in applied mathematics from the University of California, Merced. His research is on effective large scale computational methods applied 
to optimal power flow problems.
\end{IEEEbiography}

\begin{IEEEbiography}[{\includegraphics[width=1in,height=1.25in,clip,keepaspectratio]{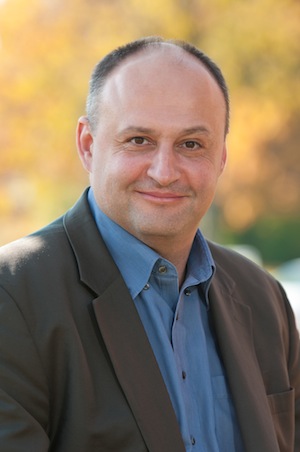}}]{Mihai Anitescu}
Is a senior computational mathematician in the Mathematics and  Computer  Science  Division  at  Argonne  National  Laboratory  and  a professor  in  the  Department  of  Statistics  at  the  University  of  Chicago.  He obtained  his  engineer  diploma  (electrical  engineering)  from  the  Polytechnic University  of  Bucharest  in  1992  and  his  Ph.D.  in  applied  mathematical  and computational sciences from the University of Iowa in 1997. He specializes in the areas of numerical optimization, computational science, numerical analysis and uncertainty quantification in which he has published more than 100 papers in  scholarly  journals  and  book  chapters.  He  is  on  the  editorial  board  of  the SIAM  Journal  on  Optimization  and  he  is  a  senior  editor  for  Optimization Methods  and  Software,  he  is  a  past  member  of  the  editorial  boards  of the  Mathematical  Programming  A  and  B,  SIAM  Journal  on  Uncertainty Quantification,  and  SIAM  Journal  on  Scientific  Computing.  He  has  been recognized for his work in applied mathematics by his selection as a SIAM Fellow in 2019.
\end{IEEEbiography}

\vspace{0.5cm}







\fbox{
\parbox{3in}{
\scriptsize
The submitted manuscript has been created by UChicago Argonne, LLC, Operator of Argonne 
National Laboratory (``Argonne''). Argonne, a U.S. Department of Energy Office of Science 
laboratory, is operated under Contract No. DE-AC02-06CH11357. The U.S. Government retains 
for itself, and others acting on its behalf, a paid-up nonexclusive, irrevocable worldwide 
license in said article to reproduce, prepare derivative works, distribute copies to the 
public, and perform publicly and display publicly, by or on behalf of the Government. 
The Department of Energy will provide public access to these results of federally 
sponsored research in accordance with the DOE Public Access 
Plan. \texttt{http://energy.gov/downloads/doe-public-accessplan}
}}
\vfill

\end{document}

%% file: cover.tex
\pagestyle{empty}
  
\vspace{1.75in}

\begin{centering}

ARGONNE NATIONAL LABORATORY

9700 South Cass Avenue

Argonne, Illinois  60439

\vspace{1.5in}

{\large \textbf{Convergence Analysis of Fixed Point Chance Constrained Optimal Power Flow Problems}}

\vspace{.5in}

\textbf{J. J. Brust, M. Anitescu}

\vspace{.5in}

Mathematics and Computer Science Division

\vspace{.25in}

Preprint ANL/MCS-P9431-0121

\vspace{.5in}

January 2022

\end{centering}

\vspace{2.0in}

\bigskip

\par\noindent
\footnotetext [1]
{
This work was supported by the U.S. Department of Energy, Office of Science,
Advanced Scientific Computing Research, under Contract DE-AC02-06CH11357
at Argonne National Laboratory.
}
\footnotetext [2]
{
J. J. Brust is now at Department of Mathematics, University of California, San Diego, CA.
}
\newpage

\vspace*{\fill}
\begin{center}
\fbox{
\parbox{4in}{
The submitted manuscript has been created by UChicago Argonne, LLC, Operator of Argonne 
National Laboratory (``Argonne''). Argonne, a U.S. Department of Energy Office of Science 
laboratory, is operated under Contract No. DE-AC02-06CH11357. The U.S. Government retains 
for itself, and others acting on its behalf, a paid-up nonexclusive, irrevocable worldwide 
license in said article to reproduce, prepare derivative works, distribute copies to the 
public, and perform publicly and display publicly, by or on behalf of the Government. 
The Department of Energy will provide public access to these results of federally 
sponsored research in accordance with the DOE Public Access 
Plan. \texttt{http://energy.gov/downloads/doe-public-accessplan}
}}
\end{center}
\vfill

\newpage
\pagestyle{plain}
\setcounter{page}{1}